\RequirePackage{fix-cm}
\documentclass[pdflatex,sn-mathphys-num]{sn-jnl}
\usepackage[T1]{fontenc}
\usepackage{lmodern}
\usepackage{amsmath,amssymb,amsfonts,amsthm,mathtools}
\usepackage{enumitem}
\usepackage{microtype}
\usepackage{xcolor}
\makeatletter
\@ifundefined{th@thmstyleone}{%
\newtheoremstyle{thmstyleone}{18pt plus2pt minus1pt}{18pt plus2pt minus1pt}%
  {\small\itshape}{0pt}{\small\bfseries}{}{.5em}%
  {\thmname{#1}\thmnumber{\@ifnotempty{#1}{ }\@upn{#2}}%
   \thmnote{ {\the\thm@notefont(#3)}}}%
\newtheoremstyle{thmstyletwo}{18pt plus2pt minus1pt}{18pt plus2pt minus1pt}%
  {\small\normalfont}{0pt}{\small\itshape}{}{.5em}%
  {\thmname{#1}\thmnumber{\@ifnotempty{#1}{ }{#2}}%
   \thmnote{ {\the\thm@notefont(#3)}}}%
\newtheoremstyle{thmstylethree}{18pt plus2pt minus1pt}{18pt plus2pt minus1pt}%
  {\small\normalfont}{0pt}{\small\bfseries}{}{.5em}%
  {\thmname{#1}\thmnumber{\@ifnotempty{#1}{ }\@upn{#2}}%
   \thmnote{ {\the\thm@notefont(#3)}}}%
}{}
\makeatother

\theoremstyle{thmstyleone}
\newtheorem{theorem}{Theorem}[section]
\newtheorem{proposition}[theorem]{Proposition}
\newtheorem{lemma}[theorem]{Lemma}
\newtheorem{corollary}[theorem]{Corollary}
\theoremstyle{thmstylethree}
\newtheorem{definition}[theorem]{Definition}

\newtheorem{question}[theorem]{Question}
\theoremstyle{thmstyletwo}
\newtheorem{remark}[theorem]{Remark}
\numberwithin{equation}{section}
\newcommand{\OO}{\mathcal O}
\newcommand{\Ok}{\mathcal O_k}
\newcommand{\Q}{\mathcal Q}
\newcommand{\PH}{P_H}
\newcommand{\PS}{P_S}
\newcommand{\Irr}{\operatorname{Irr}}
\newcommand{\sat}{\operatorname{sat}}
\newcommand{\Boxx}{\mathord{\Box}}
\newcommand{\Dia}{\mathord{\Diamond}}
\newcommand{\cl}{\operatorname{cl}}
\newcommand{\Sier}{\mathbb S}
\newcommand{\Topz}{\mathbf{Top}_0}
\newcommand{\NKV}[2]{N(#1\to #2)}
\newcommand{\doi}[1]{\url{https://doi.org/#1}}
\newcommand{\arxivlink}[1]{\href{https://arxiv.org/abs/#1}{\nolinkurl{arXiv:#1}}}
\hypersetup{hidelinks,bookmarksdepth=2,
 pdftitle={Wilker property, function spaces and hyperspaces},
 pdfauthor={Xiaoquan Xu},
 pdfsubject={Results in Mathematics submission},
 pdfkeywords={Wilker property, compact-open function space, hyperspace, Scott topology, H-sobriety, core-compactness}}
\makeatletter
\def\ps@headings{\let\@oddhead\@empty\let\@evenhead\@empty
 \def\@oddfoot{\hfil\thepage\hfil}\let\@evenfoot\@oddfoot
 \let\@mkboth\markboth}
\def\ps@titlepage{\let\@oddhead\@empty\let\@evenhead\@empty
 \def\@oddfoot{\hfil\thepage\hfil}\let\@evenfoot\@oddfoot}
\makeatother
\begin{document}
\title[Wilker property, function spaces and hyperspaces]{Wilker property, function spaces and hyperspaces}
\author*[1]{\fnm{Xiaoquan} \sur{Xu}}\email{xiqxu2002@163.com}
\affil[1]{\orgdiv{School of Computer Information Engineering},
 \orgname{Nanchang Institute of Technology},
 \orgaddress{\city{Nanchang}, \postcode{330044}, \country{China}}}
\abstract{The Wilker property requires that a compact set covered by two open sets admit compact covering pieces inside those open sets. We study this property through products, compact-open function spaces and hyperspaces. We construct countable compact Fr\'echet--Urysohn $T_1$ hereditarily Wilker spaces with non-Wilker products and squares. For spaces whose compact subsets are closed, Wilker-ness is equivalent to every compact subspace being Hausdorff; this yields a positive product theorem when one Wilker factor has closed compact subsets. For an irreducible subset system $H$, we characterize preservation of $H$-sobriety by compact-open function spaces in terms of $H$-supercompactness of compact subsets of the domain, and show that Wilker spaces satisfy this criterion. The Hoare hyperspace of a space is Wilker exactly when binary union on the Scott space of its open-set lattice is jointly continuous for the ordinary product topology. Consequently, the Hoare construction preserves and reflects Wilker-ness on consonant bases, although preservation fails in general. We also identify core-compactness of the base with local compactness of its Hoare hyperspace and prove that the Smyth construction reflects Wilker-ness.}
\maketitle
\noindent\textbf{Acknowledgements.} This research was supported by the National Natural Science Foundation of China (Nos.~12471070 and 12071199).
\par

\section{Introduction}

The compact-cover splitting condition studied in this paper has two classical sources.  It appears already in Choquet's theory of capacities, where a family of sets satisfying the corresponding binary refinement property is called \emph{rich} \cite{Choquet1954}.  It was later isolated by Wilker as condition (D) in his study of adjoint product and hom functors \cite{Wilker1970}.  We follow the established terminology and call a space $X$ \emph{Wilker} if, whenever a compact set $K$ is contained in the union of two open sets $U_1,U_2$, one can cover $K$ by compact sets $K_i\subseteq U_i$.  Wilker proved, in particular, that every Hausdorff space and every locally compact space has this property.

The condition also occurs in measure-extension theory \cite{KeimelLawson2005} and in compact-open constructions for multifunction spaces \cite{Papadopoulos1989}.  More recently, every quasi-Polish space was shown to be Wilker \cite{deBrechtKawai2019}, every LCS-complete space is Wilker \cite{dBGJL2019}, and co-consonance implies Wilker-ness through strong compactness properties \cite{ZhaoHeWang2023}.

We develop the basic theory of the Wilker property and examine its
connections with compact-open function spaces and hyperspaces. The
hyperspaces considered here are the Hoare and Smyth powerspaces,
equipped with the lower and upper Vietoris topologies, respectively.
These constructions connect structural questions in general topology
with the order-theoretic study of open-set lattices.

We first study the topology of the class itself.  Wilker spaces are hereditary with respect to open, closed and locally closed subspaces, are preserved by retracts, and are preserved and reflected by arbitrary topological sums.  Arbitrary heredity, however, fails in a particularly strong way: every $T_0$-space is homeomorphic to a subspace of a locally hypercompact Wilker space.  A new separation-theoretic feature emerges when Wilker-ness is combined with the KC property: we prove that
\[
  X\text{ is KC and Wilker}
  ~\Longleftrightarrow~
  \text{every compact subspace of }X\text{ is Hausdorff}.
\]
This characterization yields a mixed product theorem.  If $X$ is Wilker and every compact subspace of $Y$ is Hausdorff, then $X\times Y$ is Wilker.  Consequently,
\[
  X,Y\text{ Wilker and one factor KC}
  ~\Longrightarrow~
  X\times Y\text{ Wilker},
\]
and in particular the product of a Wilker space with a Hausdorff space is Wilker.  The KC hypothesis cannot be omitted. We construct two countable compact
Fr\'echet--Urysohn $T_1$ anti-Hausdorff hereditarily Wilker spaces whose
product is not Wilker (Theorem~\ref{thm:Wilker-product-negative}). A finite
sum gives a hereditarily Wilker space with a non-Wilker square; a
co-finite-core construction gives an anti-Hausdorff square counterexample
as well. Consequently, the full category of Wilker spaces is neither
reflective nor epireflective in $\Topz$.

The counterexample uses the co-finite modifications studied by Xu and
Ji in their preprint \emph{Products of co-sober spaces need not be co-sober}
(2026, \arxivlink{2610.14434v1}) for co-sobriety. Its additional ingredient is
that every proper closed subspace of the particular Kov\'ar--van Douwen
seed is compact Hausdorff (Lemma~\ref{lem:proper-closed-seed}). This
makes the singleton modifications hereditarily Wilker. An explicit
binary open cover of a compact diagonal then has no compact refinement:
the projections of any proposed pieces would express the original
anti-Hausdorff KC-space as a union of two proper closed sets. The
modifications themselves are credited to the earlier work; the
Wilker analysis and the compact-cover obstruction are established here.
Together with the mixed product theorem, this locates failure of
productivity precisely outside the guarantee supplied by a KC factor.

We next compare Wilker-ness with core-compactness.  The rational line is Wilker but not core-compact.  In the opposite direction, we show that every subspace of $B\times J$ is Wilker whenever $B$ is Hausdorff and $J=[0,1)$ carries the Scott topology of the reverse usual order.  Applied to the classical Hofmann--Lawson construction, this yields a second-countable core-compact Wilker $T_0$-space that is not locally compact.  Hence
\[
  \text{locally compact}
  ~\Longrightarrow~
  \text{core-compact and Wilker}
\]
is a strict implication, while the question whether every core-compact $T_0$-space is Wilker remains open.

The Scott-topological consequences make the connection with domain theory more explicit.  The Scott space of every quasicontinuous poset is locally compact \cite{MaoXu2006}; hence it is Wilker.  For complete lattices one can go further: since their Scott spaces are well-filtered, core-compactness already forces local compactness and therefore Wilker-ness \cite{XuShenXiZhao2020,ZhaoHeWang2023}.  On the other hand, complete lattices have coherent Scott spaces, and Zhao, He and Wang showed that coherent Wilker spaces are co-sober.  Combining this with the recent construction of Ji and Xu of a complete Boolean algebra with a non-co-sober Scott space, given in their preprint
\emph{Scott spaces of complete Boolean algebras need not be co-sober}
(2026, \arxivlink{2608.08484v1}), yields a complete Boolean algebra whose Scott space is not Wilker.  A parallel phenomenon occurs for upper topologies: the dual of Isbell's lattice, together with the lower-topology example of Wen and Xu and the coherence theorem of Xu, Yang and Chen, yields a complete lattice whose upper topology is not Wilker.  Thus strong order-completeness alone does not force Wilker-ness for either of these standard one-sided order topologies.  The unrestricted Boolean, Heyting and meet-continuous Scott-space versions all fail, while the corresponding countable questions for complete Heyting algebras and complete meet-continuous lattices remain open.

The second part of the paper turns from the basic theory to compact-open function spaces.  For a $T_0$-space $X$, let
\[
  \Ok(X)=(\OO(X),\phi(\OO(X))),
\]
where $\phi(\OO(X))$ is generated by the sets
\[
  \Phi(K)=\{U\in\OO(X):K\subseteq U\},
\]
with $K$ compact in $X$.  Given an irreducible subset system $H$, we call a subset $K\subseteq X$ \emph{$H$-supercompact} if every $H$-family in $\Ok(X)$ whose union covers $K$ already contains one member covering $K$.  For a fixed $X$, we prove the equivalence of the following four conditions:

\begin{enumerate}[label=\textup{(\arabic*)}]
\item $\Ok(X)$ is $H$-sober;
\item every compact subset of $X$ is $H$-supercompact;
\item for every $T_0$-space $Y$ the compact-open topology on $C(X,Y)$ has property $S$ with respect to $H$;
\item $C_k(X,Y)$ is $H$-sober for every $H$-sober space $Y$.
\end{enumerate}

This reorganizes several compact-open function-space questions from \cite{XuBaoZhang2023} into a single structural statement.

Wilker's condition is the main bridge to this theory.  We prove that every compact subset of a Wilker space is $H$-supercompact for every irreducible subset system $H$.  Consequently, if $X$ is Wilker, then $\Ok(X)$ is $H$-sober for every $H$, the compact-open topology on $C(X,Y)$ has property $S$ with respect to every $H$, and $C_k(X,Y)$ is $H$-sober whenever $Y$ is.  In particular, $\Ok(X)$ is sober for every Wilker space $X$.  If $X$ is also consonant, then $\phi(\OO(X))=\sigma(\OO(X))$, and the sobriety of the Scott space $\Sigma\OO(X)$ follows immediately.

The third part concerns the Hoare and Smyth hyperspace constructions.  For the nonempty Hoare powerspace we prove
\[
  X\text{ core-compact}
  \Leftrightarrow \PH(X)\text{ locally compact}
  \Leftrightarrow \PH(X)\text{ core-compact}.
\]
We further show that every $\PH(X)$ is coherent.  Since coherent Wilker spaces are co-sober, a known countable zero-dimensional metrizable space $M$ with non-co-sober $\PH(M)$ gives
\[
  M\text{ Wilker},~ \PH(M)\text{ not Wilker}.
\]
Thus the Hoare construction does not preserve Wilker-ness in general.
A further characterization identifies exactly what the Wilker property
of $\PH(X)$ measures:
\[
 \PH(X)\text{ Wilker}
 \iff
 \cup:\Sigma\OO(X)\times\Sigma\OO(X)\to\Sigma\OO(X)
       \text{ is continuous}.
\]
The domain here is the ordinary product of Scott spaces, not the Scott
topology of the order product. The proof combines the compact-family
correspondence of de Brecht and Kawai
\cite[Theorem~6.10]{deBrechtKawai2019} with a coherent-space refinement
lemma. Wilker-ness of $X$ itself is characterized by continuity of union
on $\Ok(X)$. Hence, for every consonant base,
\[
 X\text{ Wilker}\iff\PH(X)\text{ Wilker}.
\]
This gives a positive solution of the Hoare-reflection problem under
consonance; the unrestricted implication remains open. The relation
between consonance, Wilker-ness and Scott continuity of union is already
recorded in \cite[Proposition~J]{GoubaultLarrecq2026Wilker}; the additional
point here is its connection with the Hoare criterion.

For the Smyth powerspace $\PS(X)$, the known local-compactness/core-compactness characterization is complemented by the reflection theorem
\[
  \PS(X)\text{ Wilker}\Longrightarrow X\text{ Wilker}.
\]
Thus the two power constructions interact with the Wilker property and with core-compactness in markedly different ways.

The paper is organized as follows. Section~\ref{sec:prelim} fixes
notation and recalls the function-space preliminaries. Section~3
develops the basic theory, the KC--Wilker characterization and the
mixed product theorem, followed by the product and square counterexamples
and the non-reflectivity consequence. Section~4 studies the relation
with core-compactness and applications to order topologies.
Sections~5 and~6 introduce $H$-supercompactness and apply Wilker's
condition to compact-open function spaces. Section~7 treats the
Hoare construction, including the Scott-union criterion and consonant
reflection; Section~8 treats the Smyth construction.
Section~\ref{sec:questions} singles out two central open problems,
Smyth preservation and the core-compactness implication, with a brief
discussion of their related directions. The final section summarizes
the resulting structural picture.

\section{Preliminaries}\label{sec:prelim}

Throughout, spaces are $T_0$ unless explicitly stated otherwise.  For a space $X$, $\OO(X)$ denotes its lattice of open subsets, ordered by inclusion.  A subset is \emph{saturated} if it is an intersection of open sets.  We write $\Q(X)$ for the set of all nonempty compact saturated subsets of $X$.

\begin{definition}[Choquet--Wilker \cite{Choquet1954,Wilker1970}]\label{def:Wilker}
A topological space $X$ is a \emph{Wilker space} if, whenever $K\subseteq X$ is compact and $U_1,U_2\in\OO(X)$ satisfy $K\subseteq U_1\cup U_2$, there are compact subsets $K_i\subseteq U_i$ such that
\[
  K\subseteq K_1\cup K_2.
\]
\end{definition}

The compact sets $K_i$ are allowed to be empty and need not be
subsets of $K$. The original condition is binary, but the finite form
is immediate.

\begin{lemma}\label{lem:finite-wilker}
If $X$ is Wilker, $K\subseteq X$ is compact, and
\[
 K\subseteq U_1\cup\cdots\cup U_n
\]
with each $U_i$ open, then there are compact $K_i\subseteq U_i$ such that
\[
 K\subseteq K_1\cup\cdots\cup K_n.
\]
\end{lemma}

\begin{proof}
Induct on $n$, using the binary Wilker property for
$U_1\cup(U_2\cup\cdots\cup U_n)$.
\end{proof}

We shall repeatedly use the following elementary fact.

\begin{lemma}\label{lem:saturation-compact}
The saturation of a compact subset of a topological space is compact.
\end{lemma}

\begin{proof}
Let $K$ be compact and let $\{U_i:i\in I\}$ be an open cover of $\sat(K)$.  It covers $K$, so finitely many members have union $U$ containing $K$.  Since $U$ is open and contains $K$, it contains $\sat(K)$.  Hence those finitely many $U_i$ cover $\sat(K)$.
\end{proof}

Thus the Wilker condition may equivalently be stated using compact saturated sets only.

For later use it is convenient to adjoin the empty compact saturated set and write
\[
  \Q_0(X)=\Q(X)\cup\{\emptyset\},
\]
ordered by the Smyth order
\[
  K\sqsubseteq L ~\Longleftrightarrow~ L\subseteq K.
\]
Binary meets in $\Q_0(X)$ are unions.

For a compact subset $K\subseteq X$, put
\[
  \Phi(K)=\{U\in\OO(X):K\subseteq U\}.
\]
The topology on $\OO(X)$ generated by the family $\{\Phi(K): K\subseteq
X\text{ compact}\}$ is denoted by $\phi(\OO(X))$, and
\[
  \Ok(X)=(\OO(X),\phi(\OO(X))).
\]
As observed in \cite{XuBaoZhang2023},
\[
  \upsilon(\OO(X))\subseteq \phi(\OO(X))\subseteq\sigma(\OO(X)),
\]
and $\Ok(X)$ has inclusion as its specialization order. The latter can
also be seen directly, since singletons are compact.

Let $H :\Topz\to\mathbf{Set}$ be an irreducible subset system in the
sense of \cite{Xu2021}: for every $T_0$-space $Z$,
\[
  \{\{z\}:z\in Z\}\subseteq H(Z)\subseteq\Irr(Z),
\]
and continuous maps preserve $H$-subsets under direct image. A space
$Z$ is $H$-sober if every $A\in H(Z)$ has a unique generic point, i.e.,
there is a unique $z\in Z$ such that
\[
  \cl_Z A=\cl_Z\{z\}.
\]
The systems $D$, \emph{RD}, \emph{WD} and $R$ yield, respectively, $d$-spaces,
well-filtered spaces (for both \emph{RD} and \emph{WD}), and sober spaces; see
\cite{Xu2021,XuShenXiZhao2020}.

For spaces $X,Y$, let $C(X,Y)$ be the set of continuous maps from $X$ to
$Y$. The compact-open topology on $C(X,Y)$ has subbasic open sets
\[
  \NKV{K}{V}=\{f\in C(X,Y):f(K)\subseteq V\},
\]
where $K\subseteq X$ is compact and $V\in\OO(Y)$. The resulting function
space is denoted by $C_k(X,Y)$.

We use property $S$ from \cite{XuBaoZhang2023}. An order-compatible
topology $\tau$ on $C(X,Y)$ has property $S$ with respect to $H$ if, for
every $\mathcal F\in H((C(X,Y),\tau))$ and every $g\in C(X,Y)$ satisfying
\[
  g^{-1}(V)=\bigcup_{f\in\mathcal F}f^{-1}(V)
  ~(V\in\OO(Y)),
\]
one has $g\in\cl_\tau\mathcal F$.

We shall use the following result proved in \cite[Theorem 36]{XuBaoZhang2023}.

\begin{proposition}\label{prop:S-Hsober}
Let $H$ be an irreducible subset system, let $X,Y$ be $T_0$-spaces, and
let $\tau$ be a topology on $C(X,Y)$ finer than the pointwise convergence
topology. If $Y$ is $H$-sober, then $(C(X,Y),\tau)$ is $H$-sober if and
only if $\tau$ has property $S$ with respect to $H$.
\end{proposition}

We shall also use the standard homeomorphism
\begin{equation}\label{eq:Sierpinski-homeo}
  C_k(X,\Sier)\longrightarrow \Ok(X),
  ~ f\longmapsto f^{-1}(\{1\}),
\end{equation}
where $\Sier=\{0,1\}$ is the Sierpi\'nski space with open sets
$\emptyset,\{1\},\Sier$; see \cite{XuBaoZhang2023}.

\section{Basic theory of Wilker spaces}

This section develops the structural properties of the Wilker
condition used throughout the paper. After elementary permanence
properties and the KC characterization, we prove the mixed product
theorem and then give counterexamples to unrestricted productivity,
already for squares. The same examples yield non-reflectivity.

\subsection{Heredity, retracts and sums}

\begin{proposition}[Open and closed heredity]\label{prop:hereditary}
Let $X$ be Wilker.
\begin{enumerate}[label=\textup{(\arabic*)}]
\item Every open subspace of $X$ is Wilker.
\item Every closed subspace of $X$ is Wilker.
\item Consequently, every locally closed subspace of $X$ is Wilker.
\end{enumerate}
\end{proposition}

\begin{proof}
Let $A\subseteq X$.

If $A$ is open, every compact $K\subseteq A$ is compact in $X$, and every open subset of $A$ is open in $X$.  The Wilker decomposition in $X$ therefore already lies in $A$.

If $A$ is closed, write $U_i=A\cap V_i$ with $V_i$ open in $X$.  If $K\subseteq A$ is compact and $K\subseteq U_1\cup U_2$, apply Wilker-ness of $X$ to $K\subseteq V_1\cup V_2$, obtaining compact $C_i\subseteq V_i$ covering $K$.  Then $K_i=C_i\cap A$ is compact, $K_i\subseteq U_i$, and $K\subseteq K_1\cup K_2$.

The locally closed case follows by composition of the first two cases.
\end{proof}

The preceding result is essentially the strongest general hereditary statement one can expect.  In fact, arbitrary heredity fails in a very strong form.

\begin{theorem}[Failure of arbitrary heredity]\label{thm:not-hereditary}
Every $T_0$-space is homeomorphic to a subspace of a locally hypercompact Wilker space.  Consequently, the class of Wilker spaces is not hereditary with respect to arbitrary subspaces.
\end{theorem}

\begin{proof}
Let $\Sier=\{0,1\}$ be the Sierpi\'nski space and let $I$ be any set.  Identify $\Sier^I$ with $\mathcal P(I)$ by sending $f\in\Sier^I$ to
\[
 \{i\in I:f(i)=1\}.
\]
Under this identification, a basis for the product topology consists of the sets
\[
 \uparrow F=\{A\subseteq I:F\subseteq A\},
 ~ F\subseteq I\ \text{finite}.
\]
Each $\uparrow F$ is supercompact.  Indeed, if an open cover of $\uparrow F$ is given, one member contains the point $F$; since every open subset of $\Sier^I$ is an upper set for inclusion, that member contains all of $\uparrow F$.  Thus $\Sier^I$ has a basis of supercompact open sets.  In particular, it is locally hypercompact, hence locally compact and therefore Wilker by Wilker's classical theorem.

Now let $X$ be a $T_0$-space.  The canonical map
\[
 e_X:X\longrightarrow \Sier^{\OO(X)},~
 e_X(x)(U)=
 \begin{cases}
 1,&x\in U,\\
 0,&x\notin U,
 \end{cases}
\]
is a topological embedding.  Indeed, it is continuous because the inverse image under the $U$-coordinate of the open singleton $\{1\}$ is exactly $U$; it is injective because $X$ is $T_0$; and these same coordinate inverse images recover the topology of $X$.  Hence every $T_0$-space embeds into a locally hypercompact Wilker space.

Finally, non-Wilker $T_0$-spaces exist.  For example, Wilker observed that the one-point compactification $\mathbb Q^+$ of the rational line is a compact KC non-Hausdorff space and does not satisfy condition~(D) \cite{Wilker1970}.  Therefore $\mathbb Q^+$ embeds into a Wilker space while itself being non-Wilker, so arbitrary heredity fails.
\end{proof}

\begin{remark}\label{rem:sierpinski-cube}
The proof shows more than non-heredity.  The Sierpi\'nski powers $\Sier^I$, which are naturally identified with the Scott spaces of the algebraic lattices $\mathcal P(I)$, form a universal source of locally hypercompact Wilker spaces for $T_0$ embeddings.  This is one place where the basic theory of Wilker spaces meets a standard domain-theoretic representation.
\end{remark}

\begin{proposition}[Retracts]\label{prop:retract}
Every retract of a Wilker space is Wilker.
\end{proposition}

\begin{proof}
Let $e:X\to Y$ and $r:Y\to X$ be continuous with $r\circ e=\mathrm{id}_X$, and suppose $Y$ is Wilker.  Let $K\subseteq X$ be compact and $K\subseteq U_1\cup U_2$ with $U_i$ open.  Then $e(K)$ is compact and
\[
 e(K)\subseteq r^{-1}(U_1)\cup r^{-1}(U_2).
\]
Choose compact $C_i\subseteq r^{-1}(U_i)$ covering $e(K)$.  Then $K_i=r(C_i)$ is compact, lies in $U_i$, and the $K_i$ cover $K$.
\end{proof}

\begin{theorem}[Topological sums]\label{thm:sums}
For a family $(X_i)_{i\in I}$ of spaces, the topological sum $\coprod_{i\in I}X_i$ is Wilker if and only if every $X_i$ is Wilker.
\end{theorem}

\begin{proof}
Each summand is clopen, so the forward implication follows from Proposition~\ref{prop:hereditary}.

Conversely, let $X=\coprod_iX_i$ with every $X_i$ Wilker, and let $K\subseteq X$ be compact.  Since the summands form an open cover of $K$, only finitely many of them meet $K$.  On each of those finitely many components apply Wilker-ness to $K\cap X_i$ and the restrictions of the two given open sets.  Finite unions of the resulting compact pieces give the required decomposition in $X$.
\end{proof}

Let $\mathbf{Wilk}_0$ denote the full subcategory of $\mathbf{Top}_0$ consisting of Wilker spaces.  Proposition~\ref{prop:retract} and Theorem~\ref{thm:sums} show that $\mathbf{Wilk}_0$ is closed under retracts and arbitrary coproducts.

\subsection{Separation properties}

Recall that a space is \emph{KC} if every compact subset is closed.
Following Wilker \cite[p.~279]{Wilker1970}, a space $X$ is called a
\emph{$KT_4$ space} if every compact subspace $K\subseteq X$ is normal:
for any disjoint subsets $F_1,F_2$ that are closed in $K$, there are
$V_1,V_2\in\OO(K)$ such that
\[
  F_1\subseteq V_1,\qquad F_2\subseteq V_2,\qquad
  V_1\cap V_2=\emptyset.
\]
Here normality does not include the $T_1$ axiom.
Wilker proved that every $KT_4$ space satisfies condition~(D), and in particular every Hausdorff space is Wilker \cite[Theorem~8]{Wilker1970}.  The following elementary characterization shows that, after imposing the KC condition, Wilker-ness exactly forces compact subspaces to be Hausdorff.

\begin{theorem}[KC--Wilker characterization]\label{thm:KC-Wilker}
For a space $X$, the following conditions are equivalent.
\begin{enumerate}[label=\textup{(\arabic*)}]
\item $X$ is KC and Wilker.
\item Every compact subspace of $X$ is Hausdorff.
\end{enumerate}
\end{theorem}

\begin{proof} (1) $\Rightarrow$ (2): Assume first that $X$ is KC and Wilker, and let $K\subseteq X$ be compact.  Since $X$ is KC, $K$ is closed in $X$; hence $K$ is Wilker by Proposition~\ref{prop:hereditary}.  Let $x,y\in K$ be distinct.  Singletons are compact, hence closed in a KC space, so
\[
 K=(K\setminus\{y\})\cup(K\setminus\{x\})
\]
is an open cover of the compact space $K$.  By Wilker-ness there are compact subsets
\[
 A\subseteq K\setminus\{y\},~
 B\subseteq K\setminus\{x\}
\]
with $K=A\cup B$.  The sets $A$ and $B$ are compact as subsets of $X$, hence closed by the KC property.  Therefore
\[
 K\setminus B~\text{and}~ K\setminus A
\]
are disjoint open neighborhoods in $K$ of $x$ and $y$, respectively.  Thus $K$ is Hausdorff.

(2) $\Rightarrow$ (1): Suppose that every compact subspace of $X$ is Hausdorff.  We first show that $X$ is KC.  Let $A\subseteq X$ be compact.  If $x\in\cl_X A\setminus A$, then $A\cup\{x\}$ is compact and hence Hausdorff by assumption.  Since $A$ is compact, it is closed in the Hausdorff subspace $A\cup\{x\}$, contradicting $x\in\cl_X A$.  Thus every compact subset of $X$ is closed.

It remains to prove that $X$ is Wilker.  Let $K\subseteq U_1\cup U_2$ be compact, with $U_1,U_2$ open in $X$.  The subspace $K$ is compact Hausdorff, hence normal.  Applying the finite shrinking lemma in $K$ to the open cover
\[
 K\cap U_1,~ K\cap U_2,
\]
we obtain closed subsets $F_i$ of $K$ such that
\[
 F_i\subseteq K\cap U_i,~ K=F_1\cup F_2.
\]
Each $F_i$ is compact, so these sets give the required Wilker decomposition.
\end{proof}

\begin{corollary}\label{cor:compact-KC}
If $X$ is compact and KC, then $X$ is Wilker if and only if $X$ is Hausdorff.
\end{corollary}

\begin{proof}
The forward implication follows from Theorem~\ref{thm:KC-Wilker}, since $X$ itself is a compact subspace.  The reverse implication is Wilker's classical Hausdorff result.
\end{proof}

\begin{remark}\label{rem:KC-Wilker}
Thus Wilker-ness by itself is far from being a separation axiom, but in the presence of KC it has a rigid separation-theoretic meaning: it is equivalent to Hausdorffness on every compact subspace.  In particular, every compact KC non-Hausdorff space is automatically non-Wilker.
\end{remark}

\subsection{Mixed products}

The preceding characterization gives a positive product theorem under
a separation hypothesis on one factor. The counterexamples in
Subsection~\ref{subsec:negative-products} will show that this hypothesis
cannot simply be omitted.

\begin{theorem}[Mixed product theorem]\label{thm:mixed-product}
Let $X$ be a Wilker space, and suppose that every compact subspace of $Y$ is Hausdorff.  Then $X\times Y$ is Wilker.
\end{theorem}

\begin{proof}
Let $K\subseteq X\times Y$ be compact and let $U_1,U_2$ be open in $X\times Y$ with
\[
 K\subseteq U_1\cup U_2.
\]
Put $C=\pi_Y(K)$.  Then $C$ is compact and hence Hausdorff by hypothesis.  Moreover, Theorem~\ref{thm:KC-Wilker} shows that $Y$ is KC.

Fix $y\in C$ and write
\[
 K_y=\{x\in X:(x,y)\in K\},~
 U_i^y=\{x\in X:(x,y)\in U_i\}.
\]
Since $Y$ is KC, $\{y\}$ is closed, so $K_y$ is compact.  Also
\[
 K_y\subseteq U_1^y\cup U_2^y.
\]
By the Wilker property of $X$, choose compact sets
\[
 A_i(y)\subseteq U_i^y~ (i=1,2)
\]
such that $K_y\subseteq A_1(y)\cup A_2(y)$.

For each $i$, compactness of $A_i(y)$ and openness of $U_i$ give an open set $V_i(y)\subseteq X$ and an open neighborhood $W_i(y)$ of $y$ in $Y$ such that
\[
 A_i(y)\subseteq V_i(y),~
 V_i(y)\times W_i(y)\subseteq U_i.
\]
If $A_i(y)=\emptyset$, take $V_i(y)=\emptyset$ and $W_i(y)=Y$.  Put
\[
 V(y)=V_1(y)\cup V_2(y).
\]
Then $K_y\subseteq V(y)$.  Hence
\[
 F_y=K\cap\bigl((X\setminus V(y))\times Y\bigr)
\]
is compact and $y\notin\pi_Y(F_y)$.  Since $Y$ is KC, the compact set $\pi_Y(F_y)$ is closed.  Therefore
\[
 W_0(y)=Y\setminus\pi_Y(F_y)
\]
is an open neighborhood of $y$, and
\begin{equation}\label{eq:mixed-local-control}
 K\cap(X\times W_0(y))\subseteq V(y)\times W_0(y).
\end{equation}

The set
\[
 O_y=C\cap W_0(y)\cap W_1(y)\cap W_2(y)
\]
is an open neighborhood of $y$ in the compact Hausdorff space $C$.  Choose an open neighborhood $N_y$ of $y$ in $C$ such that
\[
 L_y:=\cl_C N_y\subseteq O_y.
\]
Then $L_y$ is compact; since $Y$ is KC, it is also closed in $Y$.  Consequently
\[
 K^{(y)}=K\cap(X\times L_y)
\]
is compact, and so is its projection
\[
 P_y=\pi_X(K^{(y)}).
\]
By \eqref{eq:mixed-local-control} and $L_y\subseteq W_0(y)$,
\[
 P_y\subseteq V_1(y)\cup V_2(y).
\]
Apply the Wilker property of $X$ once more to obtain compact sets
\[
 D_i(y)\subseteq V_i(y)~(i=1,2)
\]
with $P_y\subseteq D_1(y)\cup D_2(y)$.  Since $L_y\subseteq W_i(y)$,
\[
 D_i(y)\times L_y\subseteq U_i,
\]
and the two compact rectangles $D_i(y)\times L_y$ cover $K\cap(X\times L_y)$.

The family $(N_y)_{y\in C}$ covers the compact space $C$.  Choose a finite subcover
\[
 N_{y_1},\ldots,N_{y_n}.
\]
For $i=1,2$, set
\[
 H_i=\bigcup_{j=1}^n D_i(y_j)\times L_{y_j}.
\]
Each $H_i$ is compact, $H_i\subseteq U_i$, and $K\subseteq H_1\cup H_2$.  Hence $X\times Y$ is Wilker.
\end{proof}

Combining Theorems~\ref{thm:KC-Wilker} and~\ref{thm:mixed-product} gives a convenient intrinsic form.

\begin{corollary}\label{cor:KC-product}
If $X$ and $Y$ are Wilker spaces and at least one of them is KC, then $X\times Y$ is Wilker.
\end{corollary}

\begin{proof}
Suppose, for example, that $Y$ is KC.  By Theorem~\ref{thm:KC-Wilker}, every compact subspace of $Y$ is Hausdorff, so Theorem~\ref{thm:mixed-product} applies.  The other case follows by symmetry.
\end{proof}

\begin{corollary}[Hausdorff factor]\label{cor:Hausdorff-product}
If $X$ is Wilker and $Y$ is Hausdorff, then $X\times Y$ is Wilker.
\end{corollary}

\begin{proof}
Every compact subspace of a Hausdorff space is Hausdorff, so this is an immediate consequence of Theorem~\ref{thm:mixed-product}.
\end{proof}

\begin{remark}\label{rem:product-obstruction}
Consequently, if the product of two Wilker spaces fails to be Wilker,
then both factors must be non-KC, and in particular both must be
non-Hausdorff. The next subsection constructs counterexamples of
exactly this type.
\end{remark}

\begin{remark}
A naive Wallace--Tychonoff argument does not establish productivity in the
non-Hausdorff case. Given compact $K\subseteq X\times Y$ covered by
finitely many open rectangles, compactness of the projections and
Tychonoff's theorem do not by themselves give a cover of $K$ by compact
rectangles subordinate to those open rectangles. The proof of
Theorem~\ref{thm:mixed-product} succeeds because the second factor makes
compact projections closed and supplies finite shrinking inside the compact
Hausdorff space $\pi_Y(K)$. The counterexample below shows that an
unrestricted compact-cover refinement cannot hold.
\end{remark}

\subsection{Failure of product and square stability}\label{subsec:negative-products}

We now prove that the separation hypothesis in the mixed product theorem
cannot be omitted. The construction uses the co-finite modifications of
Xu and Ji \textup{(\arxivlink{2610.14434v1}, Sections~3 and~6)}, originally used to
study products of co-sober spaces. We retain those constructions and
analyze their Wilker property; the obstruction to splitting a compact
cover in the product is proved directly.

A space is \emph{hereditarily Wilker} if all of its subspaces are Wilker.
A nonempty space is \emph{anti-Hausdorff} if any two nonempty open sets
intersect. We use the particular countably infinite compact
Fr\'echet--Urysohn anti-Hausdorff KC-space $(M,\tau)$ obtained from
van Douwen's construction \cite{vanDouwen1993} in the realization described
by Kov\'ar \cite[Example~2.2 and the following paragraph, p.~6]{Kovar2005};
see also the product preprint \textup{(\arxivlink{2610.14434v1}, Lemma~5.1)}. In this realization there are
a family $\mathcal A$ of infinite subsets of $M$ and a map
$\ell:\mathcal A\to M$, with $\ell(A)\notin A$, such that a closed
subbasis for $\tau$ is
\begin{equation}\label{eq:seed-closed-subbasis}
 \mathcal S=\{B\cup\{\ell(A)\}: B\subseteq A\in\mathcal A\}.
\end{equation}
Here we use the specific realization with the stated topological
properties, not an arbitrary choice of an almost-disjoint family and
of the map $\ell$. The following additional consequence of the closed
subbasis is what is needed for Wilker-ness.

\begin{lemma}\label{lem:proper-closed-seed}
Every proper closed subspace of $M$ is compact Hausdorff. Consequently,
for each $a\in M$, every compact subspace of $M\setminus\{a\}$ is
Hausdorff.
\end{lemma}

\begin{proof}
Let $S=B\cup\{\ell(A)\}$ be a member of \eqref{eq:seed-closed-subbasis}.
For $b\in B$, the set
$(B\setminus\{b\})\cup\{\ell(A)\}$ is closed in $M$, so $\{b\}$ is open
in $S$. Thus every point of $S$ except possibly $\ell(A)$ is isolated.
Since $S$ is $T_1$, an isolated point $b$ and any other point are separated
by $\{b\}$ and $S\setminus\{b\}$. Hence $S$ is Hausdorff. It is also
compact because it is closed in the compact space $M$.

Let $F\subsetneq M$ be closed and nonempty. Choose a point of
$M\setminus F$ and a basic open neighborhood of it disjoint from $F$.
The complement of this neighborhood is a finite union
$S_1\cup\cdots\cup S_n$ of members of \eqref{eq:seed-closed-subbasis}
containing $F$. A finite union $Y$ of closed Hausdorff subspaces $S_i$
is Hausdorff: each diagonal $\Delta_{S_i}$ is closed in $S_i^2$, hence
in $Y^2$, and $\Delta_Y=\bigcup_{i=1}^n\Delta_{S_i}$ is closed.
Therefore $F$ is Hausdorff. It is compact by closedness in $M$; the
empty case is immediate.

Finally, a compact subset $K$ of $M\setminus\{a\}$ is compact in $M$,
and hence closed in $M$ by the KC property. Since it omits $a$, it is
proper, so the first assertion applies.
\end{proof}

Let more generally $(M,\tau)$ be a $T_1$-space and $a\in M$. Following the product preprint
\textup{(\arxivlink{2610.14434v1}, Section~3)}, write $\tau_{\mathrm{cof}}$ for the
co-finite topology, including the empty set, and put
\begin{equation}\label{eq:singleton-modification}
 \tau_a=\tau_{\mathrm{cof}}\cup\{U\in\tau:a\notin U\},
 \qquad X_a=(M,\tau_a).
\end{equation}
This is a compact $T_1$ topology contained in $\tau$. Every subset
containing $a$ is compact, because every neighborhood of $a$ is
co-finite. The restrictions of $\tau_a$ and $\tau$ to
$M\setminus\{a\}$ coincide: for $U\in\tau$, the set $U\setminus\{a\}$
is $\tau$-open and belongs to $\tau_a$. In particular, a subset avoiding
$a$ is compact in $X_a$ if and only if it is compact in $(M,\tau)$.
These are the singleton cases recorded in the product preprint \textup{(\arxivlink{2610.14434v1}, Lemma~3.1)}.

\begin{lemma}\label{lem:singleton-Wilker}
The space $X_a$ is Wilker if and only if the original subspace
$M\setminus\{a\}$ is Wilker. If every compact subspace of
$M\setminus\{a\}$ is Hausdorff, then $X_a$ is hereditarily Wilker.
\end{lemma}

\begin{proof}
Necessity of the first assertion follows from open heredity. Conversely,
suppose that $M\setminus\{a\}$ is Wilker, and let
$K\subseteq U_1\cup U_2$, where $K$ is compact in $X_a$ and the $U_i$
are open in $X_a$. If $a\in U_1$, put
\[
 C_1=(K\cap U_1)\cup\{a\},\qquad C_2=K\setminus U_1.
\]
The set $C_1$ is compact because it contains $a$, and $C_2$ is finite
because $U_1$ is co-finite. They lie in $U_1,U_2$, respectively, and
cover $K$. The case $a\in U_2$ is symmetric. If neither open set
contains $a$, apply Wilker-ness in the original open subspace
$M\setminus\{a\}$. Notice that the compact pieces need not be subsets
of $K$, as in Definition~\ref{def:Wilker}.

For the hereditary assertion, take $S\subseteq X_a$. If $a\notin S$,
every compact subset of $S$ is a compact Hausdorff subset of
$M\setminus\{a\}$. Theorem~\ref{thm:KC-Wilker} makes $S$ Wilker.
If $a\in S$, every relative neighborhood of $a$ is co-finite in $S$,
and every subset of $S$ containing $a$ is compact. The argument above
applies inside $S$, using the already established Wilker property of
$S\setminus\{a\}$.
\end{proof}

\begin{theorem}[Failure of productivity]\label{thm:Wilker-product-negative}
There exist two countably infinite compact Fr\'echet--Urysohn $T_1$
anti-Hausdorff hereditarily Wilker spaces whose product is not Wilker.
\end{theorem}

\begin{proof}
Take the space $(M,\tau)$ used in Lemma~\ref{lem:proper-closed-seed}
and distinct $a,b\in M$. By Lemmas~\ref{lem:proper-closed-seed}
and~\ref{lem:singleton-Wilker}, the modifications $X_a,X_b$ are
hereditarily Wilker. They are compact and $T_1$, and they remain
anti-Hausdorff because their topologies are coarser than $\tau$.
They are also Fr\'echet--Urysohn, as proved in the product preprint
\textup{(\arxivlink{2610.14434v1}, Lemma~3.3)}. Indeed, at the modified point a
sequence of distinct points of any infinite set converges to that
point; away from it, closure can be tested in the original open
complement, where the Fr\'echet--Urysohn property applies.

The diagonal map from $(M,\tau)$ to $X_a\times X_b$ is continuous.
Thus its image
\[
 D=\{(x,x):x\in M\}
\]
is compact. The two open sets
\begin{equation}\label{eq:Wilker-diagonal-cover}
 W_1=(M\setminus\{a\})\times M,\qquad
 W_2=M\times(M\setminus\{b\})
\end{equation}
cover $D$, since $a\ne b$.

Suppose that compact sets $C_i\subseteq W_i$ covered $D$.
The projection $A=\pi_1(C_1)$ is compact in $X_a$ and avoids $a$,
so it is compact in $(M,\tau)$ and therefore is a proper
$\tau$-closed set. Similarly, $B=\pi_2(C_2)$ is a proper
$\tau$-closed set avoiding $b$. Every $(x,x)\in D$ belongs to $C_1$
or $C_2$, so $M=A\cup B$. This contradicts anti-Hausdorffness of $M$.
Hence the cover \eqref{eq:Wilker-diagonal-cover} has no Wilker
refinement.
\end{proof}

The compact sets $C_i$ in the proof are allowed to contain points outside
$D$. No intersection $C_i\cap D$ is assumed compact. Also, the proof
does not infer failure of Wilker-ness merely from the known failure of
co-sobriety of the product: it excludes a compact-cover splitting
directly by taking projections.

\begin{corollary}[A non-Wilker square]\label{cor:Wilker-square-negative}
There is a countable compact Fr\'echet--Urysohn $T_1$ hereditarily
Wilker space whose square is not Wilker.
\end{corollary}

\begin{proof}
Let $Z=X_a\amalg X_b$ for the two spaces in
Theorem~\ref{thm:Wilker-product-negative}. A finite topological sum
preserves compactness, the Fr\'echet--Urysohn property and $T_1$.
Every subspace of $Z$ is a sum of subspaces of the two factors, so
Theorem~\ref{thm:sums} gives hereditary Wilker-ness. The product
$X_a\times X_b$ is a clopen subspace of $Z^2$. If $Z^2$ were Wilker,
Proposition~\ref{prop:hereditary} would contradict
Theorem~\ref{thm:Wilker-product-negative}.
\end{proof}

The square can also be made anti-Hausdorff. We use the co-finite-core
construction in the product preprint \textup{(\arxivlink{2610.14434v1}, Section~6)} and record its
additional Wilker permanence property.

\begin{proposition}\label{prop:Wilker-cofinite-core}
Let $P,Q$ be disjoint $T_1$ Wilker spaces and $R$ an infinite set
disjoint from them. On $Z=P\amalg Q\amalg R$ put
\begin{equation}\label{eq:Wilker-cofinite-core}
 \lambda=\{\emptyset\}\cup
 \{U\cup V\cup(R\setminus F):
 U\in\OO(P),\ V\in\OO(Q),\ F\subseteq R\text{ finite}\}.
\end{equation}
Then $(Z,\lambda)$ is Wilker. Consequently, a countable compact
Fr\'echet--Urysohn $T_1$ anti-Hausdorff Wilker space can have a
non-Wilker square.
\end{proposition}

\begin{proof}
The basic properties of \eqref{eq:Wilker-cofinite-core} are established
in the product preprint \textup{(\arxivlink{2610.14434v1}, Lemma~6.1)}: $P,Q$ are closed subspaces with
their original topologies, and $K\subseteq Z$ is compact if and only
if $K\cap P$ and $K\cap Q$ are compact in their respective spaces.
In particular, every subset of $R$ is compact.

Let $K\subseteq O_1\cup O_2$ be a compact set covered by two open
sets of $Z$. Split $K\cap P$ using compact covering pieces
$A_i\subseteq P\cap O_i$, and similarly split $K\cap Q$ using compact
$B_i\subseteq Q\cap O_i$. Put
\[
 E_1=K\cap R\cap O_1,\qquad
 E_2=(K\cap R)\setminus O_1\subseteq O_2.
\]
Then $A_i\cup B_i\cup E_i$ are compact subsets of $O_i$ covering $K$.

For the last assertion, take disjoint copies of $X_a,X_b$ from
Theorem~\ref{thm:Wilker-product-negative} and countable $R$.
The construction preserves compactness and the Fr\'echet--Urysohn
property and gives a $T_1$ anti-Hausdorff space
\textup{(\arxivlink{2610.14434v1}, Lemma~6.1)}. Its square contains the closed
subspace $P\times Q$, which is not Wilker. Closed heredity completes
the proof.
\end{proof}

\subsection{Categorical consequence}

\begin{theorem}\label{theor:Wilker-nonreflective}
The full subcategory $\mathbf{Wilk}_0$ of $\Topz$ is neither reflective
nor epireflective. More precisely, the square in
Corollary~\ref{cor:Wilker-square-negative} has no reflection into
Wilker spaces.
\end{theorem}

\begin{proof}
Let $Z$ be as in Corollary~\ref{cor:Wilker-square-negative}, and suppose
that $\eta:Z^2\to C$ were a reflection into Wilker spaces. Since $Z$
is Wilker, its coordinate projections factor uniquely through $\eta$:
there are continuous $p_i:C\to Z$ with $p_i\eta=\pi_i$.
Then $r=(p_1,p_2):C\to Z^2$ satisfies $r\eta=\mathrm{id}_{Z^2}$.
The uniqueness clause of the reflection, applied to $\eta$, gives
$\eta r=\mathrm{id}_C$, since both sides have composite $\eta$ with
$\eta$. Thus $\eta$ is a homeomorphism, contradicting that $C$ is
Wilker and $Z^2$ is not. Hence $\mathbf{Wilk}_0$ is not reflective,
and therefore cannot be epireflective.
\end{proof}

This is the product obstruction to reflectivity, expressed directly by
the universal property. Together with Theorem~\ref{thm:not-hereditary},
it settles the unrestricted heredity and reflection questions
negatively. The behavior under quotients and quotient-like maps is
not settled by these arguments and remains a separate direction.

\section{Sources of Wilker spaces and the core-compactness problem}

Wilker proved that every $KT_4$ space satisfies condition (D), hence every Hausdorff space is Wilker; he also proved that every locally compact space is Wilker \cite{Wilker1970}.  Later work produced large non-Hausdorff classes.  Every quasi-Polish space is Wilker \cite{deBrechtKawai2019}, every LCS-complete space is Wilker \cite{dBGJL2019}, and every co-consonant $T_0$-space is Wilker by the strong-compactness argument of \cite{ZhaoHeWang2023}.

Recall that a space $X$ is \emph{Alexandroff} if arbitrary intersections
of open subsets of $X$ are open. Equivalently, every point $x\in X$ has
a least open neighborhood, namely
\[
  \uparrow x=\bigcap\{U\in\OO(X):x\in U\},
\]
where $\uparrow x$ is taken in the specialization order. Thus the open
sets of an Alexandroff $T_0$-space are exactly the upper sets of its
specialization order.

For completeness we record two elementary sources; under our standing $T_0$ convention the zero-dimensional case is already covered by Hausdorffness, while the Alexandroff case is genuinely non-Hausdorff.

\begin{proposition}\label{prop:zeroalex}
Every zero-dimensional space is Wilker.  Every Alexandroff space is Wilker.
\end{proposition}

\begin{proof}
Let $X$ be zero-dimensional, $K\subseteq U_1\cup U_2$ compact.  For each $x\in K$ choose a clopen neighborhood $V_x$ contained in one of the $U_i$.  Compactness gives finitely many such $V_x$.  Let $C_i$ be the union of those assigned to $U_i$ and put $K_i=K\cap C_i$.  Each $K_i$ is closed in $K$, hence compact, lies in $U_i$, and the two sets cover $K$.

Now let $X$ be Alexandroff.  By Lemma~\ref{lem:saturation-compact} we may assume that $K$ is compact saturated.  Since $K$ is an upper set and the principal upper sets $\{\uparrow x:x\in K\}$ form an open cover of $K$, compactness yields a finite $F\subseteq K$ such that $K=\uparrow F$.  Assign each point of $F$ to one of the $U_i$ containing it, and let $F_i$ be the corresponding finite subsets.  Then $K_i=\uparrow F_i$ are compact, lie in $U_i$, and cover $K$.
\end{proof}

\subsection{Core-compactness}

Recall that $X$ is core-compact if $\OO(X)$ is a continuous lattice.  Every locally compact space is core-compact, and for sober spaces core-compactness is equivalent to local compactness \cite{GHKLS2003}.  Hence local compactness gives
\[
  \text{locally compact}\Longrightarrow
  \text{core-compact and Wilker}.
\]
Neither condition alone characterizes local compactness; the examples below also distinguish their conjunction from local compactness.

The Wilker property does not imply core-compactness: the rational line $\mathbb Q$ is Hausdorff, hence Wilker, but it is not locally compact; being Hausdorff, it is sober, so it cannot be core-compact.

The opposite implication is one of the central open problems of this paper.

\begin{question}[Core-compact versus Wilker]\label{q:CCW}
Is every core-compact $T_0$-space Wilker?
\end{question}

The classical Hofmann--Lawson example does not provide a negative answer.  In fact it yields a rather strong positive test case.  We first prove a general lemma.

Let $J=[0,1)$ equipped with the Scott topology of the reverse usual order.  Thus the open sets of $J$ are $[0,a)$, $0\leq a\leq 1$.

\begin{theorem}\label{thm:BxJ}
Let $B$ be a Hausdorff space.  Every subspace $X$ of $B\times J$ is Wilker.
\end{theorem}

\begin{proof}
Let $K\subseteq X$ be compact and let $U_1,U_2$ be open in $X$ with $K\subseteq U_1\cup U_2$.  Write $U_i=X\cap O_i$ with $O_i$ open in $B\times J$.  Let $\pi:B\times J\to B$ be projection and put $C=\pi(K)$.  Then $C$ is compact Hausdorff.

For $x\in C$, the fiber
\[
 K_x=\{t\in J:(x,t)\in K\}
\]
is nonempty and compact: indeed, $K_x$ is homeomorphic to the closed subspace $K\cap(\{x\}\times J)$ of $K$, because $B$ is Hausdorff.  Every nonempty compact subset of $J$ has a largest element in the usual order.  Indeed, if $s=\sup K_x$ were not attained (or were equal to $1$), the increasing family of open initial intervals $[0,a)$ with $a<s$ (respectively $a<1$) would cover $K_x$ without a finite subcover.  Write
\[
 q(x)=\max K_x.
\]
The map $q:C\to[0,1)$ is upper semicontinuous.  For $r\in[0,1)$,
\[
 \{x\in C:q(x)\ge r\}
 =\pi\bigl(K\cap(B\times[r,1))\bigr).
\]
Here $[r,1)$ is closed in $J$, since its complement $[0,r)$ is Scott open.  Hence the set inside the projection is compact; its image is compact in the Hausdorff space $B$, and therefore closed.  Thus $q$ is upper semicontinuous in the ordinary real-valued sense.

For $i=1,2$ and $x\in C$, let
\[
\begin{aligned}
\mathcal A_i(x)=\{a\in[0,1]:{}&\ \exists\,N\in\mathcal N_B(x)\text{ such that}\\
&X\cap(N\times[0,a))\subseteq U_i\},\\
f_i(x)={}&\sup\mathcal A_i(x).
\end{aligned}
\]
where $\mathcal N_B(x)$ denotes the open neighborhoods of $x$ in $B$.
Each $f_i$ is lower semicontinuous: if $f_i(x)>r$, one may choose $a>r$ and a witnessing neighborhood $N$; then $N\cap C\subseteq\{f_i>r\}$.

Put
\[
 W_i=\{x\in C:q(x)<f_i(x)\}.
\]
Since $q$ is upper semicontinuous and $f_i$ is lower semicontinuous, $W_i$ is open in $C$.  Moreover $C=W_1\cup W_2$.  Indeed, for $x\in C$ the point $(x,q(x))$ belongs to $K$, hence to one of the $U_i$.  By openness, there are $a>q(x)$ and an open neighborhood $N$ of $x$ such that $X\cap(N\times[0,a))\subseteq U_i$, so $q(x)<f_i(x)$.

Since $C$ is compact Hausdorff, the finite shrinking lemma gives closed sets $F_i\subseteq W_i$ with $C=F_1\cup F_2$.  Define
\[
 K_i=K\cap\pi^{-1}(F_i).
\]
Then each $K_i$ is compact and $K=K_1\cup K_2$.  We claim $K_i\subseteq U_i$.  If $(x,t)\in K_i$, then $t\le q(x)<f_i(x)$.  Choose $c$ with $q(x)<c<f_i(x)$.  By the definition of the supremum, some $a>c$ is witnessed by a neighborhood $N$ of $x$ with $X\cap(N\times[0,a))\subseteq U_i$.  Since $t<a$, we obtain $(x,t)\in U_i$.
\end{proof}

Hofmann and Lawson constructed a second-countable core-compact $T_0$-space $X_{HL}$ which is not locally compact; in fact every compact subset of $X_{HL}$ has empty interior \cite{HofmannLawson1978,GHKLS2003}.  Their space is a subspace of $[0,1]\times J$.

\begin{corollary}\label{cor:HL}
There exists a second-countable core-compact Wilker $T_0$-space that is not locally compact.
\end{corollary}

\begin{proof}
Apply Theorem~\ref{thm:BxJ} to the Hofmann--Lawson space.
\end{proof}

Thus
\[
 \{\text{locally compact }T_0\text{-spaces}\}
 \subsetneq
 \{\text{core-compact Wilker }T_0\text{-spaces}\}.
\]
In particular, even the conjunction of core-compactness and Wilker-ness does not force local compactness.

\subsection{Scott spaces and domain-theoretic consequences}\label{subsec:Scott-Wilker}

We now record several consequences for Scott spaces.  They provide a direct link between the basic theory developed above and standard classes from domain theory.

Recall that Mao and Xu call a poset $P$ \emph{quasicontinuous} if, for every $x\in U\in\sigma(P)$, there is a nonempty finite set $F\subseteq P$ such that
\[
 x\in\operatorname{int}_{\sigma(P)}\uparrow F
 \subseteq \uparrow F\subseteq U.
\]
It is clear that the Scott space of every quasicontinuous poset is locally compact (cf. \cite[Proposition~2.6]{MaoXu2006}).  Wilker's classical theorem therefore gives the following immediate consequence.

\begin{corollary}\label{cor:quasicont-Wilker}
If $P$ is a quasicontinuous poset, then its Scott space $\Sigma P$ is Wilker.
\end{corollary}

In particular, the conclusion applies to quasicontinuous domains, continuous posets and continuous domains.  Notice that directed completeness is not required in Corollary~\ref{cor:quasicont-Wilker}.

The core-compact Scott-space problem has a complete answer for complete lattices.

\begin{proposition}\label{prop:CC-complete-lattice-Wilker}
Let $L$ be a complete lattice.  If $\Sigma L$ is core-compact, then $\Sigma L$ is locally compact, sober and Wilker.
\end{proposition}

\begin{proof}
The Scott space of every complete lattice is well-filtered; see, for example, \cite{Xi-Lawson-2017}.  Every core-compact well-filtered $T_0$-space is locally compact (and sober) by \cite{XuShenXiZhao2020}.  Since every locally compact space is Wilker, the conclusion follows.
\end{proof}

Thus the specialization of Question~\ref{q:CCW} to complete lattices is positive.  For general Scott spaces, however, a natural domain-theoretic version of the question remains.

\begin{question}[Core-compact Scott spaces]\label{q:Scott-corecompact}
Let $P$ be a poset, or even a dcpo.  If its Scott space $\Sigma P$ is core-compact, must $\Sigma P$ be Wilker?
\end{question}

We next record a useful obstruction.  Recall that a $T_0$-space is \emph{co-sober} if every $k$-irreducible compact saturated subset is the saturation of a point.  Zhao, He and Wang proved that every Wilker $T_0$-space is weakly co-sober and that, in a coherent space, weak co-sobriety implies co-sobriety \cite{ZhaoHeWang2023}.  Hence:

\begin{proposition}\label{prop:coherent-Wilker-cosober}
Every coherent Wilker $T_0$-space is co-sober.
\end{proposition}

Since the Scott space of every complete lattice is coherent, we immediately obtain the following complete-lattice form.

\begin{corollary}\label{prop:complete-lattice-Wilker-cosober}
Let $L$ be a complete lattice.  If $\Sigma L$ is Wilker, then $\Sigma L$ is co-sober.
\end{corollary}

This turns recent non-co-sobriety examples into Wilker counterexamples without any additional analysis of compact covers.  Ji and Xu constructed a complete Boolean algebra $B$ such that $\Sigma B$ is not co-sober; in fact their example is also non-sober and well-filtered \textup{(\arxivlink{2608.08484v1})}.  Corollary~\ref{prop:complete-lattice-Wilker-cosober} therefore yields:

\begin{corollary}\label{cor:Boolean-nonWilker}
There exists a complete Boolean algebra whose Scott space is not Wilker.
Consequently, there exist a complete Heyting algebra and a complete meet-continuous lattice whose Scott spaces are not Wilker.
\end{corollary}

\begin{remark}\label{rem:distributive-nonWilker}
The counterexample in Corollary~\ref{cor:Boolean-nonWilker} is much stronger than a merely distributive complete-lattice example: every complete Boolean algebra is a frame, hence a complete Heyting algebra, and in particular a complete meet-continuous lattice.  Thus none of these three unrestricted classes has Wilker Scott spaces automatically.
\end{remark}

Countability changes the picture at the Boolean level.  Ji and Xu proved that every countable complete Boolean algebra is finite \textup{(\arxivlink{2608.08484v1}, Proposition~5.4)}.  Hence its Scott topology is finite Alexandroff and therefore Wilker by Proposition~\ref{prop:zeroalex}.  The next two cases are substantially less clear.

\begin{question}[Countable Scott--Wilker problems]\label{q:countable-Scott-Wilker}
\leavevmode
\begin{enumerate}[label=\textup{(\arabic*)}]
\item Is the Scott space of every countable complete Heyting algebra (equivalently, every countable frame) Wilker?
\item Is the Scott space of every countable complete meet-continuous lattice Wilker?
\end{enumerate}
\end{question}

The second problem admits an even broader variant in which ``complete meet-continuous lattice'' is replaced by ``countable meet-continuous dcpo''.  These countable questions are not consequences of Corollary~\ref{cor:Boolean-nonWilker}, since the Boolean case collapses to the finite case.  They therefore isolate a genuinely countable boundary between the positive quasicontinuous theory and the unrestricted Boolean counterexample.

\subsection{Upper topologies of complete lattices}\label{subsec:upper-Wilker}

There is a parallel negative phenomenon for the upper topology.  For a poset $P$, let $\nu(P)$ denote the upper topology generated by the subbasic sets
\[
  P\setminus\downarrow x,~ x\in P.
\]
Equivalently, $\nu(P)=\omega(P^{\mathrm{op}})$, where $\omega$ denotes the lower topology.

Wen and Xu proved that Isbell's complete lattice $I$ equipped with its lower topology $\omega(I)$ is sober but not co-sober \cite{WenXu2018}.  Let $L=I^{\mathrm{op}}$.  Then $L$ is again a complete lattice and
\[
  (L,\nu(L))=(I,\omega(I))
\]
as topological spaces.  On the other hand, Xu, Yang and Chen proved that the upper topology of every sup-complete poset is coherent (indeed, it satisfies considerably stronger properties) \cite{XuYangChen2026}.  Hence $(L,\nu(L))$ is coherent.  Proposition~\ref{prop:coherent-Wilker-cosober} now gives the following.

\begin{corollary}\label{cor:upper-nonWilker}
There exists a complete lattice whose upper topology is not Wilker.
\end{corollary}

\begin{proof}
For the lattice $L=I^{\mathrm{op}}$ above, $(L,\nu(L))$ is not co-sober but is coherent.  If it were Wilker, Proposition~\ref{prop:coherent-Wilker-cosober} would make it co-sober, a contradiction.
\end{proof}

\begin{remark}
Corollaries~\ref{cor:Boolean-nonWilker} and \ref{cor:upper-nonWilker} give two complementary complete-lattice obstructions:
\[
\begin{array}{c|c}
\text{Scott topology} & \text{even complete Boolean algebras may fail to be Wilker},\\
\text{upper topology} & \text{even complete lattices may fail to be Wilker}.
\end{array}
\]
Thus Wilker-ness is not forced merely by completeness, even for standard order-compatible topologies.
\end{remark}

\section{\texorpdfstring{$H$}{H}-supercompactness and compact-open function spaces}\label{sec:Hsuper}

We now introduce the graded compactness notion that controls the compact-open topology.

\begin{definition}\label{def:Hsupercompact}
Let $H$ be an irreducible subset system and let $X$ be a $T_0$-space. A
subset $K\subseteq X$ is called \emph{$H$-supercompact} if, for every
$\mathcal A\in H(\Ok(X))$,
\[
  K\subseteq\bigcup\mathcal A
  ~\Longrightarrow~
  K\subseteq U\text{ for some }U\in\mathcal A.
\]
When $H=R$, we also say that $K$ is \emph{irreducibly supercompact}.
\end{definition}

\begin{remark}\label{rem:monotone}
If $H_1\leq H_2$, then every $H_2$-supercompact subset is an
$H_1$-supercompact subset. Hence, for the standard systems,
\[
\begin{aligned}
  R\text{-supercompact}
  &\Longrightarrow WD\text{-supercompact}
   \Longrightarrow RD\text{-supercompact}\\
  &\Longrightarrow D\text{-supercompact}.
\end{aligned}
\]
Ordinary supercompactness implies $R$-supercompactness.
\end{remark}

The directed level is exactly ordinary compactness.

\begin{proposition}\label{prop:Dcompact}
For a subset $K$ of a $T_0$-space $X$, the following are equivalent:
\begin{enumerate}[label=\textup{(\arabic*)}]
\item $K$ is compact;
\item $K$ is $D$-supercompact.
\end{enumerate}
\end{proposition}

\begin{proof}
Suppose first that $K$ is compact, and let $\mathcal D\in D(\Ok(X))$
with $K\subseteq\bigcup\mathcal D$. Then a finite subfamily
$U_1,\dots,U_n\in\mathcal D$ covers $K$. Since $\mathcal D$ is directed
under inclusion, there exists $U\in\mathcal D$ with
$U_1\cup\cdots\cup U_n\subseteq U$. Hence $K\subseteq U$.

Conversely, assume that $K$ is $D$-supercompact and let
$\mathcal U\subseteq\OO(X)$ be an open cover of $K$. Let $\mathcal D$
be the directed family of all finite unions of members of $\mathcal U$.
Then $K\subseteq\bigcup\mathcal D$, so $K\subseteq V$ for some
$V\in\mathcal D$. Thus finitely many members of $\mathcal U$ cover $K$.
\end{proof}

The next observation is the technical link between compact-open function
spaces and $\Ok(X)$.

\begin{lemma}\label{lem:preimage-map}
Let $X,Y$ be $T_0$-spaces and let $V\in\OO(Y)$. Then the map
\[
  \theta_V:C_k(X,Y)\longrightarrow\Ok(X),
  ~ \theta_V(f)=f^{-1}(V),
\]
is continuous.
\end{lemma}

\begin{proof}
For every compact $K\subseteq X$,
\[
  \theta_V^{-1}(\Phi(K))
  =\{f\in C(X,Y):K\subseteq f^{-1}(V)\}
  =\NKV{K}{V},
\]
which is open in the compact-open topology.
\end{proof}

We now obtain the central equivalence.

\begin{theorem}\label{thm:main-equivalence}
Let $H$ be an irreducible subset system and let $X$ be a $T_0$-space.
Then the following conditions are equivalent:
\begin{enumerate}[label=\textup{(\arabic*)}]
\item $\Ok(X)$ is $H$-sober;
\item every compact subset of $X$ is $H$-supercompact;
\item for every $T_0$-space $Y$, the compact-open topology on $C(X,Y)$
      has property $S$ with respect to $H$;
\item for every $H$-sober space $Y$, the function space $C_k(X,Y)$ is
      $H$-sober.
\end{enumerate}
\end{theorem}

\begin{proof}
$(1)\Rightarrow(2)$. Let $K\subseteq X$ be compact and let
$\mathcal A\in H(\Ok(X))$ satisfy
$K\subseteq\bigcup\mathcal A$. By the $H$-sobriety of $\Ok(X)$, there
is a unique $U_{\mathcal A}\in\OO(X)$ such that
\[
  \cl_{\Ok(X)}\mathcal A
  =\cl_{\Ok(X)}\{U_{\mathcal A}\}.
\]
Since the specialization order of $\Ok(X)$ is inclusion, every
$U\in\mathcal A$ satisfies $U\subseteq U_{\mathcal A}$, whence
$\bigcup\mathcal A\subseteq U_{\mathcal A}$.

For the reverse inclusion, let $x\in U_{\mathcal A}$. Then
$U_{\mathcal A}\in\Phi(\{x\})$. Since
$U_{\mathcal A}\in\cl_{\Ok(X)}\mathcal A$, the open set
$\Phi(\{x\})$ meets $\mathcal A$. Hence $x\in U$ for some
$U\in\mathcal A$. Thus
\[
  U_{\mathcal A}=\bigcup\mathcal A.
\]
Now $U_{\mathcal A}\in\Phi(K)$, and again
$U_{\mathcal A}\in\cl_{\Ok(X)}\mathcal A$ yields
$\Phi(K)\cap\mathcal A\neq\emptyset$. Therefore
$K\subseteq U$ for some $U\in\mathcal A$.

$(2)\Rightarrow(1)$. Let $\mathcal A\in H(\Ok(X))$ and put
$U_{\mathcal A}=\bigcup\mathcal A$. Since every
$U\in\mathcal A$ is contained in $U_{\mathcal A}$,
\[
  \cl_{\Ok(X)}\mathcal A
  \subseteq \cl_{\Ok(X)}\{U_{\mathcal A}\}.
\]
We show that $U_{\mathcal A}\in\cl_{\Ok(X)}\mathcal A$. Let
$\Phi(K)$ be a basic neighbourhood of $U_{\mathcal A}$, where $K$ is
compact. Then
\[
  K\subseteq U_{\mathcal A}=\bigcup\mathcal A.
\]
By (2), $K\subseteq U$ for some $U\in\mathcal A$, so
$U\in\Phi(K)$. Hence every basic neighbourhood of $U_{\mathcal A}$
meets $\mathcal A$, and therefore
\[
  \cl_{\Ok(X)}\mathcal A
  =\cl_{\Ok(X)}\{U_{\mathcal A}\}.
\]
The generic point is unique because $\Ok(X)$ is $T_0$.

$(2)\Rightarrow(3)$. Let $Y$ be a $T_0$-space,
$\mathcal F\in H(C_k(X,Y))$, and $g\in C(X,Y)$ satisfy
\[
  g^{-1}(V)=\bigcup_{f\in\mathcal F}f^{-1}(V)
  ~(V\in\OO(Y)).
\]
Let $\NKV{K}{V}$ be a subbasic compact-open neighbourhood of $g$.
By Lemma~\ref{lem:preimage-map} and functoriality of $H$,
\[
  \mathcal A_V:=\{f^{-1}(V):f\in\mathcal F\}
  \in H(\Ok(X)).
\]
Moreover,
\[
  K\subseteq g^{-1}(V)=\bigcup\mathcal A_V.
\]
By (2), there is $f\in\mathcal F$ such that
$K\subseteq f^{-1}(V)$, equivalently
$f\in\NKV{K}{V}$. Thus every subbasic neighbourhood of $g$ meets
$\mathcal F$. Since $\mathcal F$ is irreducible, every finite
intersection of such neighbourhoods also meets $\mathcal F$. Hence
$g\in\cl_{C_k(X,Y)}\mathcal F$, proving property $S$.

$(3)\Rightarrow(4)$ follows from Proposition~\ref{prop:S-Hsober}, since
the compact-open topology is finer than the pointwise convergence
topology.

$(4)\Rightarrow(1)$. The Sierpi\'nski space $\Sier$ is sober, hence
$H$-sober. Therefore $C_k(X,\Sier)$ is $H$-sober by (4). By the
homeomorphism \eqref{eq:Sierpinski-homeo}, $\Ok(X)$ is $H$-sober.
\end{proof}

The theorem immediately reorganizes several of the questions from
\cite{XuBaoZhang2023}.

\begin{corollary}\label{cor:old-questions}
For a fixed $T_0$-space $X$ the following hold.
\begin{enumerate}[label=\textup{(\arabic*)}]
\item The following are equivalent:
  \begin{enumerate}[label=\textup{(\alph*)}]
  \item $\Ok(X)$ is sober;
  \item every compact subset of $X$ is $R$-supercompact;
  \item for every $T_0$-space $Y$, the compact-open topology on
        $C(X,Y)$ has property $S$ with respect to $R$;
  \item for every sober space $Y$, $C_k(X,Y)$ is sober.
  \end{enumerate}
\item The following are equivalent:
  \begin{enumerate}[label=\textup{(\alph*)}]
  \item $\Ok(X)$ is well-filtered;
  \item every compact subset of $X$ is \emph{WD}-supercompact;
  \item every compact subset of $X$ is \emph{RD}-supercompact;
  \item for every $T_0$-space $Y$, the compact-open topology on
        $C(X,Y)$ has property $S$ with respect to \emph{WD} (equivalently,
        with respect to \emph{RD});
  \item for every well-filtered space $Y$, $C_k(X,Y)$ is well-filtered.
  \end{enumerate}
\end{enumerate}
Consequently, the universal forms of Questions 49, 55, 58 and 60 in
\cite{XuBaoZhang2023} are equivalent, and the universal forms of
Questions 48, 56, 57 and 59 are equivalent (for either \emph{WD} or \emph{RD}).
\end{corollary}

\begin{proof}
Apply Theorem~\ref{thm:main-equivalence} with $H=R$, \emph{WD} and \emph{RD},
and use the characterizations of sober and well-filtered spaces by
$R$-, \emph{WD}- and \emph{RD}-sobriety from \cite{Xu2021,XuShenXiZhao2020}.
\end{proof}

\begin{corollary}\label{cor:WD-RD}
For every $T_0$-space $X$, the following are equivalent:
\begin{enumerate}[label=\textup{(\arabic*)}]
\item every compact subset of $X$ is \emph{WD}-supercompact;
\item every compact subset of $X$ is \emph{RD}-supercompact.
\end{enumerate}
\end{corollary}

\begin{remark}
Corollary~\ref{cor:WD-RD} concerns the assertion that \emph{all} compact
subsets have the corresponding property. It does not assert that the
\emph{WD} and \emph{RD} notions of supercompactness coincide for an
individual fixed subset.
\end{remark}

There is a close relation with the \emph{SI}-compactness introduced by Zhao
and Ho \cite{ZhaoHo2015}. They call $X$ $SI$-compact if every open cover
of $X$ which is irreducible in the Scott space $\Sigma\OO(X)$ contains
$X$ itself. Since
$\phi(\OO(X))\subseteq\sigma(\OO(X))$, every subset irreducible in the
Scott topology is irreducible in $\phi(\OO(X))$. Hence
$R$-supercompactness in the sense of Definition~\ref{def:Hsupercompact}
is at least as strong as the corresponding Scott-irreducible covering
condition. If $X$ is consonant, then
$\phi(\OO(X))=\sigma(\OO(X))$, and the two irreducibility tests agree.

\section{Wilker property and \texorpdfstring{$H$}{H}-supercompactness}\label{sec:Wilker-Hsuper}

We now show that Wilker's condition is a strong bridge from ordinary compactness to the whole $H$-supercompactness hierarchy. The Wilker definition and its finite form were given in Section~\ref{sec:prelim}; we therefore pass directly to the main implication.

\begin{theorem}\label{thm:Wilker-Hsuper}
Let $X$ be a Wilker space and let $H$ be an irreducible subset system.
Then every compact subset of $X$ is $H$-supercompact.
\end{theorem}

\begin{proof}
Let $K\subseteq X$ be compact and let
$\mathcal A\in H(\Ok(X))$ satisfy
$K\subseteq\bigcup\mathcal A$. Since $K$ is compact, there are
$U_1,\dots,U_n\in\mathcal A$ such that
\[
  K\subseteq U_1\cup\cdots\cup U_n.
\]
By Lemma~\ref{lem:finite-wilker}, choose compact
$K_i\subseteq U_i$ such that
\[
  K\subseteq K_1\cup\cdots\cup K_n.
\]
For each $i$, the basic open set $\Phi(K_i)$ of $\Ok(X)$ meets
$\mathcal A$, because $U_i\in\mathcal A\cap\Phi(K_i)$. Since
$\mathcal A\in H(\Ok(X))\subseteq\Irr(\Ok(X))$, irreducibility gives
\[
  \mathcal A\cap\bigcap_{i=1}^n\Phi(K_i)\neq\emptyset.
\]
Choose $U$ in this intersection. Then $K_i\subseteq U$ for every $i$,
and therefore $K\subseteq U$. Thus $K$ is $H$-supercompact.
\end{proof}

Combining Theorems~\ref{thm:main-equivalence} and
\ref{thm:Wilker-Hsuper} yields the main function-space consequence.

\begin{corollary}\label{cor:Wilker-main}
Let $X$ be a Wilker space and let $H$ be an irreducible subset system.
Then we have the following conclusions:
\begin{enumerate}[label=\textup{(\arabic*)}]
\item $\Ok(X)$ is $H$-sober.
\item For every $T_0$-space $Y$, the compact-open topology on $C(X,Y)$
      has property $S$ with respect to $H$.
\item For every $H$-sober space $Y$, $C_k(X,Y)$ is $H$-sober.
\end{enumerate}
\end{corollary}

\begin{corollary}\label{cor:Wilker-special}
Let $X$ be a Wilker space. Then we have the following conclusions:
\begin{enumerate}[label=\textup{(\arabic*)}]
\item $\Ok(X)$ is sober (hence well-filtered).
\item If $Y$ is sober, then $C_k(X,Y)$ is sober.
\item If $Y$ is well-filtered, then $C_k(X,Y)$ is well-filtered.
\item The compact-open topology on $C(X,Y)$ has property $S$ with
      respect to $R$, \emph{WD} and \emph{RD}, for every $T_0$-space $Y$.
\end{enumerate}
In particular, Questions 48, 49, 55, 56, 59 and 60 of
\cite{XuBaoZhang2023} have affirmative answers whenever the domain
space $X$ is Wilker.
\end{corollary}

Wilker showed that condition (D) (i.e., the condition of Definition~\ref{def:Wilker}) holds for some important classes of spaces,
including Hausdorff spaces and locally compact spaces
\cite{Wilker1970}. Thus Corollary~\ref{cor:Wilker-main} applies to all
such domain spaces.

We also obtain a short proof of a recent result on Scott-open set
lattices.

\begin{corollary}\label{cor:consonant-Wilker}
If $X$ is consonant and Wilker, then the Scott space
$\Sigma\OO(X)$ is sober.
\end{corollary}

\begin{proof}
By Corollary~\ref{cor:Wilker-main}, $\Ok(X)$ is sober. Consonance is
equivalent to
\[
  \phi(\OO(X))=\sigma(\OO(X)),
\]
see, for example, \cite{GoubaultLarrecq2013}.
Hence $\Sigma\OO(X)=\Ok(X)$ is sober.
\end{proof}

\begin{remark}
Corollary~\ref{cor:consonant-Wilker} recovers the consonant-Wilker part
of Z.~He's preprint \emph{When do weakly first-countable spaces and the Scott topology of open set lattice become sober?}
(2025, \arxivlink{2503.12766v2}). The present proof separates the two mechanisms:
Wilker's condition supplies the irreducible compact-cover refinement,
while consonance identifies the compact-generated topology
$\phi(\OO(X))$ with the Scott topology.
\end{remark}

\section{The Hoare power construction}

\subsection{Coherence, local compactness and non-preservation}

Let $\PH(X)$ denote the set of all nonempty closed subsets of $X$ equipped with the lower Vietoris topology, generated by
\[
 \Dia U=\{F\in\PH(X):F\cap U\ne\emptyset\},
 ~ U\in\OO(X).
\]
Schalk proved that the Hoare hyperspace (with the empty closed set included) is always sober \cite{Schalk1993}; the nonempty part is the open subspace $\Dia X$, hence it is sober as well.

\begin{proposition}[Schalk]\label{prop:PH-sober}
For every $T_0$-space $X$, $\PH(X)$ is sober.
\end{proposition}

The lower Vietoris semilattice also carries a useful compactness property that is independent of the base space.

\begin{proposition}\label{prop:PH-coherent}
For every $T_0$-space $X$, the Hoare powerspace $\PH(X)$ is coherent.
\end{proposition}

\begin{proof}
Put $Y=\PH(X)$ and consider the union map
\[
  m:Y\times Y\longrightarrow Y,~ m(A,B)=A\cup B.
\]
For every $U\in\OO(X)$,
\[
  m^{-1}(\Dia U)
  =(\Dia U\times Y)\cup(Y\times\Dia U),
\]
so $m$ is continuous.  Let $\mathcal K_1,\mathcal K_2$ be compact saturated subsets of $Y$.  Since the specialization order of $Y$ is inclusion, both families are upward closed under inclusion.  Therefore
\[
  m(\mathcal K_1\times\mathcal K_2)=\mathcal K_1\cap\mathcal K_2.
\]
Indeed, $A\cup B$ belongs to both families whenever $A\in\mathcal K_1$ and $B\in\mathcal K_2$, while every $C\in\mathcal K_1\cap\mathcal K_2$ equals $m(C,C)$.  The intersection is thus a continuous image of a compact space, hence compact; it is saturated as an intersection of saturated sets.  Thus $Y$ is coherent.
\end{proof}

Combining Proposition~\ref{prop:PH-coherent} with Proposition~\ref{prop:coherent-Wilker-cosober} yields a useful obstruction.

\begin{corollary}\label{cor:PH-Wilker-cosober}
For every $T_0$-space $X$,
\[
  \PH(X)\text{ Wilker}\Longrightarrow\PH(X)\text{ co-sober}.
\]
\end{corollary}

For later use we include a direct compact-neighborhood argument in the locally compact case.  Schalk proved the corresponding lower-Vietoris local-compactness result \cite[Proposition~6.11]{Schalk1993}; stronger formulations for core-compact bases are also standard \cite[Proposition~2.9]{HofmannNora2015}.

\begin{lemma}\label{lem:PH-LC}
If $X$ is locally compact, then $\PH(X)$ is locally compact.
\end{lemma}

\begin{proof}
Let $F\in\PH(X)$ and let
\[
 \mathcal U=\Dia U_1\cap\cdots\cap\Dia U_n
\]
be a basic neighborhood of $F$.  Choose $x_i\in F\cap U_i$.  By local compactness, choose open $V_i$ and compact saturated $Q_i$ such that
\[
 x_i\in V_i\subseteq Q_i\subseteq U_i.
\]
Then
\[
 F\in\bigcap_i\Dia V_i
 \subseteq \bigcap_i\Dia Q_i
 \subseteq \mathcal U.
\]
It remains to show that $\mathcal K=\bigcap_i\Dia Q_i$ is compact.

Consider
\[
 \varphi:Q_1\times\cdots\times Q_n\longrightarrow\PH(X),
 ~
 \varphi(q_1,\ldots,q_n)=\overline{\{q_1,\ldots,q_n\}}.
\]
For open $U\subseteq X$,
\[
 \varphi^{-1}(\Dia U)=\bigcup_{i=1}^n\pi_i^{-1}(U),
\]
because an open set meets $\overline{\{q_1,\ldots,q_n\}}$ if and only if it contains at least one $q_i$.  Hence $\varphi$ is continuous.  Its image is compact.  Since the specialization order of $\PH(X)$ is inclusion, the saturation of this image consists exactly of the closed sets containing $\overline{\{q_1,\ldots,q_n\}}$ for some $(q_1,\ldots,q_n)\in Q_1\times\cdots\times Q_n$; equivalently, it is $\bigcap_i\Dia Q_i=\mathcal K$.  By Lemma~\ref{lem:saturation-compact}, $\mathcal K$ is compact.
\end{proof}

We now relate the Hoare construction to core-compactness.  Let $X^s$ be the sobrification of $X$.  Closed subsets of $X$ and $X^s$ correspond naturally: for closed $F\subseteq X$, put
\[
 \widehat F=\{A\in\Irr_c(X):A\subseteq F\}.
\]

\begin{lemma}\label{lem:PH-sobrification}
The map $F\mapsto\widehat F$ is a homeomorphism
\[
 \PH(X)\cong\PH(X^s).
\]
\end{lemma}

\begin{proof}
The standard sobrification induces an order isomorphism $\OO(X)\cong\OO(X^s)$ and hence, by complementation, a bijection between the closed subsets of $X$ and those of $X^s$; under this bijection $F$ is sent to $\widehat F$.  It remains only to compare the lower Vietoris topologies.  If $U\in\OO(X)$ corresponds to
\[
 U^s=\{A\in\Irr_c(X):A\cap U\ne\emptyset\},
\]
then
\[
 \widehat F\cap U^s\ne\emptyset
 ~\Longleftrightarrow~
 F\cap U\ne\emptyset.
\]
Indeed, the forward implication is immediate.  Conversely, if $x\in F\cap U$, then $\overline{\{x\}}\subseteq F$ is irreducible closed and meets $U$.  Thus lower Vietoris subbasic opens correspond exactly.
\end{proof}

The implication from core-compactness of $X$ to local compactness of the full lower Vietoris hyperspace is contained in the classical theory recalled above.  The following formulation for the nonempty Hoare powerspace also records the converse in a form useful below.

\begin{theorem}\label{thm:PH-CC-LC}
For every $T_0$-space $X$, the following are equivalent:
\begin{enumerate}[label=\textup{(\arabic*)}]
\item $X$ is core-compact;
\item $\PH(X)$ is locally compact;
\item $\PH(X)$ is core-compact.
\end{enumerate}
\end{theorem}

\begin{proof}
$(1)\Rightarrow(2)$.  The sobrification $X^s$ of a core-compact space is locally compact \cite{GHKLS2003}.  By Lemma~\ref{lem:PH-LC}, $\PH(X^s)$ is locally compact, and by Lemma~\ref{lem:PH-sobrification} it is homeomorphic to $\PH(X)$.

$(2)\Rightarrow(1)$.  Let $x\in U\in\OO(X)$.  The canonical map
\[
 \eta_X:X\to\PH(X),~ x\mapsto\overline{\{x\}},
\]
is continuous and satisfies $\eta_X^{-1}(\Dia U)=U$.  By local compactness of $\PH(X)$, choose an open $\mathcal V$ and a compact set $\mathcal K$ such that
\[
 \eta_X(x)\in\mathcal V\subseteq\mathcal K\subseteq\Dia U.
\]
Put $V=\eta_X^{-1}(\mathcal V)$.  We show $V\ll U$ in $\OO(X)$.  Let $(U_i)_{i\in I}$ be directed with $U\subseteq\bigcup_iU_i$.  Then
\[
 \mathcal K\subseteq\Dia U\subseteq\bigcup_i\Dia U_i.
\]
Compactness and directedness yield $i_0$ with $\mathcal K\subseteq\Dia U_{i_0}$.  Hence
\[
 V\subseteq\eta_X^{-1}(\mathcal K)\subseteq U_{i_0}.
\]
Thus every $x\in U$ has a neighborhood $V\ll U$, so $X$ is core-compact.

$(2)\Leftrightarrow(3)$.  Local compactness always implies core-compactness.  Conversely, $\PH(X)$ is sober by Proposition~\ref{prop:PH-sober}, and a sober core-compact space is locally compact \cite{GHKLS2003}.
\end{proof}

\begin{corollary}\label{cor:PH-CC}
The Hoare powerspace construction preserves and reflects core-compactness.
\end{corollary}

\begin{corollary}\label{cor:PH-LC}
The Hoare powerspace construction preserves local compactness but does not reflect local compactness.
\end{corollary}

\begin{proof}
Preservation is Lemma~\ref{lem:PH-LC}.  For non-reflection, take the Hofmann--Lawson space $X_{HL}$.  It is not locally compact, but it is core-compact, so $\PH(X_{HL})$ is locally compact by Theorem~\ref{thm:PH-CC-LC}.
\end{proof}

The reason for this asymmetry is transparent: $\PH(X)$ carries sobriety automatically.  The Hoare construction therefore turns the core-compactness of $X$ into core-compactness plus sobriety of $\PH(X)$, which is exactly local compactness.

Combining Theorem~\ref{thm:PH-CC-LC} with Wilker's classical theorem gives the following large positive region.

\begin{corollary}\label{cor:CC-PHW}
If $X$ is core-compact, then $\PH(X)$ is locally compact, sober and Wilker.
\end{corollary}

This places Question~\ref{q:CCW} in a new light: for every core-compact $X$, the Hoare powerspace is already Wilker; the missing issue is whether the Wilker property can be recovered on the base space.

The preservation problem, however, has a strong negative answer. In their preprint \emph{A countable metric space whose Hoare power space is not co-sober}
(2026, \arxivlink{2610.05868v1}), Xu and Ji constructed a countable zero-dimensional metrizable space $M$ such that $\PH(M)$ is not co-sober.  Since $M$ is metrizable, it is Hausdorff and hence Wilker, whereas Corollary~\ref{cor:PH-Wilker-cosober} excludes Wilker-ness of its Hoare powerspace.

\begin{corollary}[Failure of Hoare preservation]\label{cor:PH-nonpreserve}
There exists a countable zero-dimensional metrizable Wilker space $M$ such that $\PH(M)$ is not Wilker.  Consequently, the Hoare powerspace construction does not preserve Wilker-ness, even within countable zero-dimensional metrizable spaces.
\end{corollary}

\begin{proof}
Take $M$ as above.  If $\PH(M)$ were Wilker, then it would be co-sober by Corollary~\ref{cor:PH-Wilker-cosober}, contrary to the metric-space preprint \textup{(\arxivlink{2610.05868v1})}.
\end{proof}

The unrestricted reflection direction remains open. The results below
give an operation-theoretic reformulation and settle it on consonant
bases.

\begin{question}[Hoare reflection]\label{q:PH-reflect}
If $\PH(X)$ is Wilker, must $X$ be Wilker?
\end{question}

A negative answer to Question~\ref{q:CCW} would automatically give a negative answer to Question~\ref{q:PH-reflect}: if $X$ were core-compact but non-Wilker, then $\PH(X)$ would be locally compact and hence Wilker by Corollary~\ref{cor:CC-PHW}.

Although preservation fails in general, it holds on important subclasses.  Core-compact bases are covered by Corollary~\ref{cor:CC-PHW}; moreover, lower powerspaces of quasi-Polish spaces are quasi-Polish, and quasi-Polish spaces are Wilker \cite{deBrechtKawai2019}.

\subsection{Union on the compact-open lattice}

We first record an elementary operation-theoretic form of the Wilker
condition. In the proof, empty compact saturated sets are allowed.

\begin{proposition}\label{prop:compact-open-union}
A space $X$ is Wilker if and only if
\[
 \cup:\Ok(X)\times\Ok(X)\longrightarrow\Ok(X)
\]
is continuous, where the domain has the ordinary product topology.
\end{proposition}

\begin{proof}
Suppose $X$ is Wilker and $U\cup V\in\Phi(K)$ for a compact saturated
$K\subseteq X$. Choose compact saturated $K_1\subseteq U$ and
$K_2\subseteq V$ with $K\subseteq K_1\cup K_2$. Then
$\Phi(K_1)\times\Phi(K_2)$ is a neighborhood of $(U,V)$ mapped into
$\Phi(K)$.

Conversely, let $K\subseteq U\cup V$ be compact saturated and suppose
union is continuous. The sets $\Phi(C)$ form a base for $\Ok(X)$,
since finite intersections correspond to finite unions of compact
sets. There are therefore compact saturated $K_1\subseteq U$ and
$K_2\subseteq V$ such that
\[
 \cup\bigl(\Phi(K_1)\times\Phi(K_2)\bigr)\subseteq\Phi(K).
\]
If $x\in K\setminus(K_1\cup K_2)$, saturation provides open sets
$U_i'$ containing $K_i$ and omitting $x$. Their union must contain
$K$, a contradiction. Thus $K\subseteq K_1\cup K_2$, proving the
Wilker condition by Lemma~\ref{lem:saturation-compact}.
\end{proof}

\subsection{A Scott-lattice criterion for Hoare powerspaces}

We next show that the corresponding operation on the Scott space
$\Sigma\OO(X)$ characterizes Wilker-ness of the Hoare powerspace.
One technical step is needed to pass from lower Vietoris subbasic
covers to arbitrary open covers.

\begin{lemma}[A coherent-space refinement lemma]\label{lem:coherent-subbasis}
Let $Y$ be coherent, and let $\mathcal S$ be an open subbasis closed
under finite unions and containing $\emptyset,Y$. Suppose that
whenever a compact saturated $K$ satisfies $K\subseteq S\cup T$ with
$S,T\in\mathcal S$, there are compact saturated sets
$A\subseteq S$, $B\subseteq T$ covering $K$. Then $Y$ is Wilker.
\end{lemma}

\begin{proof}
By Lemma~\ref{lem:saturation-compact}, it suffices to split compact
saturated sets. We may replace the compact pieces supplied by the
hypothesis with their intersections with $K$, using coherence.

Let $K\subseteq O_1\cup O_2$ with $K$ compact saturated and $O_i$
open. Choose a finite cover of $K$ by basic opens, each contained
in one $O_i$, and group them into two finite unions. Distributivity
and finite-union closure of $\mathcal S$ express these unions as
\[
 G_1=\bigcap_{j=1}^m S_j\subseteq O_1,\qquad
 G_2=\bigcap_{k=1}^n T_k\subseteq O_2,\qquad
 K\subseteq G_1\cup G_2,
\]
where $S_j,T_k\in\mathcal S$. An empty group means that $K$ is
contained in the other open set and is immediate; a whole-space group
is equally immediate. We may therefore take $m,n\geq1$.

For each $j,k$, use $K\subseteq S_j\cup T_k$ to obtain compact
saturated $A_{jk},B_{jk}\subseteq K$ covering $K$, with
$A_{jk}\subseteq S_j$ and $B_{jk}\subseteq T_k$. Put
\[
 C_1=\bigcap_{j=1}^m\bigcup_{k=1}^n A_{jk},\qquad
 C_2=\bigcup_{j=1}^m\bigcap_{k=1}^n B_{jk}.
\]
Finite unions and coherence make both sets compact saturated, and
$C_1\subseteq G_1$, $C_2\subseteq G_2$. If $x\in K\setminus C_1$,
then some $j$ satisfies $x\notin A_{jk}$ for every $k$, so
$x\in B_{jk}$ for every $k$ and hence $x\in C_2$.
Thus $K\subseteq C_1\cup C_2$.
\end{proof}

For the next correspondence we temporarily include the empty closed
set. Let $\mathcal A(X)$ be the space of all closed subsets of $X$
with the lower Vietoris topology, and use $\Dia U$ for
$\{F\in\mathcal A(X):F\cap U\ne\emptyset\}$ when working in
$\mathcal A(X)$. Then $\PH(X)=\Dia X$ is its open subspace of
nonempty closed sets.

\begin{lemma}\label{lem:Hoare-empty-point}
The space $\mathcal A(X)$ is Wilker if and only if $\PH(X)$ is
Wilker. Moreover, $\mathcal A(X)$ is coherent.
\end{lemma}

\begin{proof}
The forward implication follows from open heredity. An open set in
$\mathcal A(X)$ containing the point $\emptyset$ is the whole
space. If a two-set cover contains such an open set, take the
entire compact set as one piece and the empty set as the other.
Otherwise the open sets and the compact set being covered lie in
$\PH(X)$, so Wilker-ness there proves the converse.

For coherence, the proof of Proposition~\ref{prop:PH-coherent}
applies unchanged with the empty closed set included: union is
continuous, and its image on the product of two compact saturated
families is their intersection.
\end{proof}

Let $\mathcal K$ be a compact saturated subset of $\mathcal A(X)$,
possibly empty, and let $\mathcal H\in\sigma(\OO(X))$. Define
\begin{equation}\label{eq:Hoare-correspondence}
 \begin{split}
 \mathfrak h(\mathcal K)
   &=\{U\in\OO(X):\mathcal K\subseteq\Dia U\},\\
 \mathfrak k(\mathcal H)
   &=\bigcap_{U\in\mathcal H}\Dia U.
 \end{split}
\end{equation}
We use the following established order correspondence, without any
countable-basis or consonance hypothesis.

\begin{proposition}[de Brecht--Kawai
{\cite[Definition~6.7 and Theorem~6.10]{deBrechtKawai2019}}]
\label{prop:Hoare-order-correspondence}
The maps in \eqref{eq:Hoare-correspondence} are mutually inverse,
inclusion-reversing bijections between the compact saturated subsets
of $\mathcal A(X)$ and the Scott-open subsets of $\OO(X)$. In
particular,
\begin{align}
 U\in\mathfrak h(\mathcal K)
   &\iff\mathcal K\subseteq\Dia U,
       \label{eq:Hoare-membership}\\
 \mathfrak h(\mathcal K_1\cup\mathcal K_2)
   &=\mathfrak h(\mathcal K_1)\cap\mathfrak h(\mathcal K_2),
       \label{eq:Hoare-union-intersection}\\
 \mathfrak k(\mathcal H_1\cap\mathcal H_2)
   &=\mathfrak k(\mathcal H_1)\cup\mathfrak k(\mathcal H_2).
       \label{eq:Hoare-intersection-union}
\end{align}
\end{proposition}

In the Smyth order on compact families these maps are order-preserving;
``inclusion-reversing'' refers to ordinary set inclusion on both sides.
The empty cases are included: $\mathfrak h(\emptyset)=\OO(X)$ and
$\mathfrak k(\emptyset)=\mathcal A(X)$. Only the order correspondence
is used below, not the stronger topological assertion in
\cite[Theorem~6.11]{deBrechtKawai2019}, which requires a countable base.

\begin{theorem}[Hoare--Scott union criterion]\label{thm:Hoare-Scott-union}
For every $T_0$-space $X$, the following are equivalent:
\begin{enumerate}[label=\textup{(\arabic*)}]
\item $\PH(X)$ is Wilker;
\item the union operation
\[
 \cup:\Sigma\OO(X)\times\Sigma\OO(X)
       \longrightarrow\Sigma\OO(X)
\]
is jointly continuous, where the domain is the ordinary product
of the two Scott spaces.
\end{enumerate}
\end{theorem}

\begin{proof}
Assume (1), and let $U\cup V\in\mathcal H\in\sigma(\OO(X))$.
By Proposition~\ref{prop:Hoare-order-correspondence}, the compact
saturated family $\mathfrak k(\mathcal H)$ satisfies
\[
 \mathfrak k(\mathcal H)\subseteq\Dia(U\cup V)
                         =\Dia U\cup\Dia V.
\]
Use Lemma~\ref{lem:Hoare-empty-point} and Wilker-ness to choose
compact saturated families $\mathcal K_1\subseteq\Dia U$ and
$\mathcal K_2\subseteq\Dia V$ covering $\mathfrak k(\mathcal H)$.
Put $\mathcal H_i=\mathfrak h(\mathcal K_i)$. These are Scott-open
neighborhoods of $U,V$, respectively, and
\eqref{eq:Hoare-union-intersection} gives
\[
 \mathcal H_1\cap\mathcal H_2\subseteq\mathcal H.
\]
For $U'\in\mathcal H_1$ and $V'\in\mathcal H_2$, upperness implies
$U'\cup V'\in\mathcal H_1\cap\mathcal H_2\subseteq\mathcal H$.
This proves joint continuity.

Conversely, assume (2). Suppose that a compact saturated family
$\mathcal K\subseteq\mathcal A(X)$ is covered by $\Dia U$ and
$\Dia V$. Then $U\cup V\in\mathfrak h(\mathcal K)$.
Continuity gives Scott-open neighborhoods $\mathcal H_1,\mathcal H_2$
of $U,V$ such that $U'\cup V'\in\mathfrak h(\mathcal K)$ whenever
$U'\in\mathcal H_1$ and $V'\in\mathcal H_2$. Taking $U'=V'$ shows
\[
 \mathcal H_1\cap\mathcal H_2\subseteq\mathfrak h(\mathcal K).
\]
Apply \eqref{eq:Hoare-intersection-union} to obtain
\[
 \mathcal K\subseteq
 \mathfrak k(\mathcal H_1)\cup\mathfrak k(\mathcal H_2),\qquad
 \mathfrak k(\mathcal H_1)\subseteq\Dia U,\quad
 \mathfrak k(\mathcal H_2)\subseteq\Dia V.
\]
These are compact saturated pieces. The diamonds, together with
$\mathcal A(X)$, form an open subbasis closed under finite unions;
a cover using $\mathcal A(X)$ itself splits trivially.
Lemma~\ref{lem:Hoare-empty-point} gives coherence, and
Lemma~\ref{lem:coherent-subbasis} upgrades this subbasic splitting
property to Wilker-ness of $\mathcal A(X)$. Another application of
Lemma~\ref{lem:Hoare-empty-point} proves (1).
\end{proof}

\begin{remark}\label{rem:union-Scott-product-distinction}
Union always preserves directed suprema as a map from the order product
$\OO(X)\times\OO(X)$. This does not imply the joint continuity in
Theorem~\ref{thm:Hoare-Scott-union}: its domain has topology
$\sigma(\OO(X))\times\sigma(\OO(X))$, not the possibly finer topology
$\sigma(\OO(X)\times\OO(X))$. Thus the theorem is not a consequence
of order-theoretic Scott continuity alone.
\end{remark}

\subsection{Consonant bases and the reflection question}

Recall that consonance is the equality
$\phi(\OO(X))=\sigma(\OO(X))$. The preceding two operation-theoretic
characterizations therefore agree on consonant base spaces.

\begin{corollary}[Hoare preservation and reflection on consonant bases]
\label{cor:consonant-Hoare-Wilker}
For every consonant $T_0$-space $X$,
\[
 X\text{ is Wilker}\quad\Longleftrightarrow\quad\PH(X)\text{ is Wilker}.
\]
In particular, Question~\ref{q:PH-reflect} has an affirmative answer
for consonant $X$, and every counterexample in the general case
must be non-consonant.
\end{corollary}

\begin{proof}
Consonance identifies $\Ok(X)$ with $\Sigma\OO(X)$. Apply
Proposition~\ref{prop:compact-open-union} and
Theorem~\ref{thm:Hoare-Scott-union}.
\end{proof}

For consonant spaces, the equivalence between Wilker-ness and joint
Scott continuity of union was already obtained by de Brecht and is
recorded in \cite[Proposition~J]{GoubaultLarrecq2026Wilker}, where it is
attributed to a conversation with Keimel. We use this established
connection; the Hoare consequence here comes from combining it with
Theorem~\ref{thm:Hoare-Scott-union}.

\begin{remark}[Compact transversals]\label{rem:Hoare-compact-transversal}
There is also a direct proof of the reflection direction that isolates
the use of consonance. Let $K\subseteq U_1\cup U_2$ be compact in
$X$, and put $\eta(x)=\cl_X\{x\}$. If $\PH(X)$ is Wilker, split
$\eta(K)$ into compact families $\mathcal K_i\subseteq\Dia U_i$.
Saturate the families, which keeps them compact and inside $\Dia U_i$.
Consonance gives compact saturated $C_i\subseteq U_i$ such that
\[
 \Phi(C_i)\subseteq\mathfrak h(\mathcal K_i).
\]
Every $F\in\mathcal K_i$ meets $C_i$: otherwise $X\setminus F$
belongs to $\Phi(C_i)$, but cannot belong to
$\mathfrak h(\mathcal K_i)$. If $\eta(x)\in\mathcal K_i$, then
$\cl_X\{x\}\cap C_i\ne\emptyset$, and saturation of $C_i$ implies
$x\in C_i$. Thus $K\subseteq C_1\cup C_2$. This is the compact
transversal mechanism of
\cite[Lemma~6.12 and Theorem~6.13]{deBrechtKawai2019}.
\end{remark}

For arbitrary $X$, Question~\ref{q:PH-reflect} is now equivalent to
whether the following implication always holds:
\begin{equation}\label{eq:Hoare-reflection-transfer}
 \begin{gathered}
 \cup:\Sigma\OO(X)\times\Sigma\OO(X)\to\Sigma\OO(X)
       \text{ is continuous}\\
 \Longrightarrow\\
 \cup:\Ok(X)\times\Ok(X)\to\Ok(X)
       \text{ is continuous}.
 \end{gathered}
\end{equation}
Both domains here have the ordinary product topology. The inclusion
$\phi(\OO(X))\subseteq\sigma(\OO(X))$ alone does not justify this
transfer of joint continuity. We have proved it when $X$ is consonant,
but do not claim it for arbitrary $X$. In particular, the general
Hoare-reflection problem and the core-compactness question
\ref{q:CCW} remain unresolved here.

\section{The Smyth power construction}

Let $\PS(X)$ be $\Q(X)$ equipped with the upper Vietoris topology, whose basic open sets are
\[
 \Boxx U=\{K\in\Q(X):K\subseteq U\},
 ~ U\in\OO(X).
\]
The specialization order is reverse inclusion.

The local-compactness behavior of the Smyth construction is very different from that of the Hoare construction.

\begin{theorem}[Lyu--Chen--Jia \cite{LyuChenJia2022}]\label{thm:PS-LC}
For every topological space $X$, the following are equivalent:
\begin{enumerate}[label=\textup{(\arabic*)}]
\item $X$ is locally compact;
\item $\PS(X)$ is locally compact;
\item $\PS(X)$ is core-compact.
\end{enumerate}
\end{theorem}

Hence local compactness gives a first positive answer to the Wilker-preservation problem for Smyth powerspaces.

\begin{corollary}\label{cor:LC-PSW}
If $X$ is locally compact, then $\PS(X)$ is Wilker.
\end{corollary}

We now prove a reflection theorem.  We use the following standard compact-union lemma.

\begin{lemma}[Schalk \cite{Schalk1993}]\label{lem:union-compact-family}
If $\mathcal K$ is a compact subset of $\PS(X)$, then
\[
 \bigcup\mathcal K
\]
is a compact saturated subset of $X$.
\end{lemma}

\begin{proof}
The union is saturated because every member of $\mathcal K$ is saturated.  Let $\{U_j:j\in J\}$ be an open cover of $\bigcup\mathcal K$.  For each $Q\in\mathcal K$, compactness of $Q$ gives a finite $J_Q\subseteq J$ such that $Q\subseteq V_Q:=\bigcup_{j\in J_Q}U_j$.  The sets $\Boxx V_Q$ form an open cover of $\mathcal K$.  Choose finitely many $Q_1,\ldots,Q_m$ whose corresponding boxes cover $\mathcal K$.  The finitely many $U_j$ appearing in the $J_{Q_k}$ cover every member of $\mathcal K$, hence cover $\bigcup\mathcal K$.
\end{proof}

\begin{theorem}[Smyth reflection of Wilker-ness]\label{thm:PS-reflect}
If $\PS(X)$ is Wilker, then $X$ is Wilker.
\end{theorem}

\begin{proof}
Let $K\subseteq X$ be compact and $K\subseteq U_1\cup U_2$ with $U_i$ open.  Consider the canonical map
\[
 \xi_X:X\to\PS(X),~ \xi_X(x)=\uparrow x.
\]
The set $\uparrow x$ is compact saturated by Lemma~\ref{lem:saturation-compact}, applied to the compact singleton $\{x\}$.  The map $\xi_X$ is continuous, since $\xi_X^{-1}(\Boxx U)=U$.  Therefore $\xi_X(K)$ is compact in $\PS(X)$, and
\[
 \xi_X(K)\subseteq \Boxx U_1\cup\Boxx U_2.
\]
Since $\PS(X)$ is Wilker, there are compact subsets $\mathcal K_i\subseteq\Boxx U_i$ such that
\[
 \xi_X(K)\subseteq\mathcal K_1\cup\mathcal K_2.
\]
Set $K_i=\bigcup\mathcal K_i$.  By Lemma~\ref{lem:union-compact-family}, $K_i$ is compact saturated, and $K_i\subseteq U_i$.  If $x\in K$, then $\uparrow x=\xi_X(x)$ belongs to one of the $\mathcal K_i$, hence $x\in K_i$.  Thus $K\subseteq K_1\cup K_2$.
\end{proof}

So the difficult direction is preservation.

\begin{question}[Smyth preservation]\label{q:PS-preserve}
If $X$ is Wilker, must $\PS(X)$ be Wilker?
\end{question}

Again there are large positive subclasses.  Besides the locally compact case from Corollary~\ref{cor:LC-PSW}, the Smyth powerspace of a quasi-Polish space is quasi-Polish, hence Wilker \cite{deBrechtKawai2019}.

The Hofmann--Lawson space provides a particularly natural test case: it is core-compact, Wilker and not locally compact.  By Theorem~\ref{thm:PS-LC}, its Smyth powerspace is not core-compact and not locally compact.  Thus Question~\ref{q:PS-preserve} is genuinely outside the automatic locally compact regime already for this classical example.

\section{Two central open problems and further directions}\label{sec:questions}

\subsection{Two central open problems on Wilker spaces}

Theorem~\ref{thm:Wilker-product-negative} settles the binary product
problem negatively, even for countable compact hereditarily Wilker
factors. Corollary~\ref{cor:Wilker-square-negative} gives a square
counterexample, and Corollary~\ref{theor:Wilker-nonreflective} settles
categorical reflectivity negatively. We therefore single out two
remaining structural problems.

\begin{description}[leftmargin=2.9em,labelwidth=2.4em,itemsep=0.9em]
\item[\textnormal{(W1)}] \textbf{Smyth preservation.}
If $X$ is Wilker, must $\PS(X)$ be Wilker?

This is Question~\ref{q:PS-preserve}. The reverse implication is
Theorem~\ref{thm:PS-reflect}, and preservation holds for locally
compact and quasi-Polish spaces. The Hofmann--Lawson space remains
a natural test case: it is core-compact and Wilker but not locally
compact, and its Smyth powerspace lies outside the automatic locally
compact regime.

\item[\textnormal{(W2)}] \textbf{Core-compactness versus Wilker-ness.}
Is every core-compact $T_0$-space Wilker?

This is Question~\ref{q:CCW}. The complete-lattice Scott-space case
is positive by Proposition~\ref{prop:CC-complete-lattice-Wilker}.
Corollary~\ref{cor:HL} shows that core-compactness together with
Wilker-ness does not force local compactness, so an affirmative answer
cannot be reduced to a general local-compactness assertion.
\end{description}

The product counterexamples do not answer either question. They give
failure of a binary permanence property, whereas (W1) concerns a
specific powerspace construction and (W2) concerns the existence of
compact-cover refinements in all core-compact spaces.

\subsection{The Hoare connection and related directions}

Question~\ref{q:PH-reflect} is now equivalent to the transfer of joint
continuity in \eqref{eq:Hoare-reflection-transfer}. It has a positive
answer on consonant bases by
Corollary~\ref{cor:consonant-Hoare-Wilker}, but its general form
remains open. This question is closely tied to (W2): a negative answer
to (W2) would give a core-compact non-Wilker space $X$ whose Hoare
powerspace is locally compact and Wilker by
Corollary~\ref{cor:CC-PHW}. Conversely, a positive answer to general
Hoare reflection would imply a positive answer to (W2).
The new criterion does not by itself supply this transfer outside
consonance.

The order-topological variants remain useful test problems rather than
additional central questions. Question~\ref{q:Scott-corecompact}
asks for (W2) on Scott spaces of arbitrary posets or dcpos;
Question~\ref{q:countable-Scott-Wilker} concerns countable frames and
countable complete meet-continuous lattices. The unrestricted Boolean
and upper-topology counterexamples do not decide those countable cases.
The quotient behavior of Wilker spaces also remains separate from
the product and reflection results proved here.

On the function-space side, it remains natural to ask whether sobriety
of $\Ok(X)$, equivalently $R$-supercompactness of every compact subset
of $X$, forces Wilker-ness. A negative answer would isolate a weakening
of compact-cover splitting. Other directions concern strictness and
possible collapses in the supercompactness hierarchy of
Remark~\ref{rem:monotone}, and detection of Wilker-type properties by
property $S$ or $H$-sobriety of $C_k(X,Y_0)$ for a fixed finite $Y_0$,
such as $\Sier^2$ or a small Alexandroff space. These connect naturally
with consonance, co-consonance and Hofmann--Mislove correspondences;
(W1) and (W2) remain the two principal open problems emphasized here.

\section{Conclusions and outlook}

Wilker's compact-cover condition admits descriptions in terms of both
topological structure and continuity of union. The KC characterization supplies
a sufficient separation hypothesis for the mixed product theorem.
The co-finite modifications, by contrast, give hereditarily Wilker
factors with a non-Wilker product and square. Thus the class is stable
under retracts and sums, but not under products, and it does not form
a reflective full subcategory of $\Topz$.

The compact-open results use $H$-supercompactness to describe the role of the Wilker property. A compact refinement on the domain yields
$H$-sobriety of $\Ok(X)$, property $S$ for compact-open function
spaces, and preservation of $H$-sobriety by $C_k(X,-)$. The order-topological
examples delimit this mechanism: quasicontinuous Scott spaces are
Wilker, whereas complete Boolean algebras need not have Wilker Scott
spaces and complete lattices need not have Wilker upper topologies.

The Hoare--Scott union criterion supplies another link between the
topology of a space and the topology of its open-set lattice.
Wilker-ness of $X$ is continuity of union on $\Ok(X)$, whereas
Wilker-ness of $\PH(X)$ is continuity of union on $\Sigma\OO(X)$.
Consonance identifies these topologies and gives preservation and
reflection on that class of bases. Outside consonance, preservation
fails and the general reflection implication is not settled here.
Smyth powerspaces always reflect Wilker-ness, but their preservation
problem remains open.

The resulting picture leaves (W1), Smyth preservation, and (W2), the
core-compactness implication, as the two main remaining questions.
The counterexamples remove productivity from that list, while the
Hoare criterion gives a precise additional route toward (W2) without
claiming an answer to it.

\clearpage
\section*{Statements and Declarations}
\bmhead{Funding}
This research was supported by the National Natural Science Foundation of China (Nos.~12471070 and 12071199).
\bmhead{Competing interests}
The author has no relevant financial or non-financial interests to disclose.
\bmhead{Author contributions}
Xiaoquan Xu is the sole author and is responsible for the research and the preparation of the manuscript.
\bmhead{Data availability}
Data sharing is not applicable to this article as no datasets were generated or analysed during the current study.
\bmhead{Use of AI-assisted tools}
During the preparation of this manuscript, the author used ChatGPT
(OpenAI) to assist with language editing, organization of the exposition,
bibliographic checks, and the development and checking of mathematical
arguments. The author is responsible for verifying all mathematical
statements, proofs and references and takes full responsibility for
the final text.
\end{document}